\documentclass[notitlepage, 11pt]{article}

\usepackage{enumitem}

\usepackage{mathtools,bbm}
\usepackage{color}
\usepackage{latexsym}
\usepackage{amsmath, amssymb, amsfonts}
\usepackage[toc]{appendix}
\numberwithin{equation}{section}

\parskip  \baselineskip
\newtheorem{lemma}{Lemma}
\newtheorem{theorem}{Theorem}

\newtheorem{example}{Example}
\newtheorem{proposition}{Proposition}

\newtheorem{assumption}{Assumption}

\newtheorem{remark}{Remark}
\newtheorem{construction}{Construction}

\newcommand{\beginsec}{
\setcounter{lemma}{0}
\setcounter{theorem}{0}
\setcounter{corollary}{0}
\setcounter{definition}{0}
\setcounter{example}{0}
\setcounter{proposition}{0}
\setcounter{condition}{0}
\setcounter{assumption}{0}
\setcounter{conjecture}{0}
\setcounter{problem}{0}
\setcounter{remark}{0}
\setcounter{construction}{0}}

\newcommand{\noi}{\noindent}

\newcommand{\E}{\mathbb{E}}
\newcommand{\R}{\mathbb{R}}

\newcommand{\N}{\mathbb{N}}

\newcommand{\la}{\lambda}

\newcommand{\eps}{\varepsilon}

\newcommand{\al}{\alpha}

\newcommand{\PP}{{\mathbb P}}

\newcommand{\calC}{{\cal C}}
\newcommand{\calD}{{\cal D}}

\newcommand{\calF}{{\cal F}}

\newcommand{\calL}{{\cal L}}
\newcommand{\calM}{{\cal M}}

\newcommand{\calP}{{\cal P}}
\newcommand{\calQ}{{\cal Q}}

\newcommand{\skp}{\vspace{\baselineskip}}

\newcommand{\To}{\Rightarrow}

\newcommand\iy{\infty}

\newcommand{\qed}{\hfill $\Box$}

\newcommand{\limn}{\lim_{n\to\iy}}

\newcommand{\1}{\boldsymbol{\mathbbm{1}}}

\title{Analysis of a Finite State Many Player Game Using its Master Equation\thanks{This is the final version of the paper. To appear in {\it SIAM Journal on Control and Optimization}}.}
\author{Erhan Bayraktar\thanks{Department of Mathematics,
University of Michigan, Ann Arbor, MI 48109, USA. Emails: erhan@umich.edu, web: http://www.math.lsa.umich.edu/$\sim$erhan/. The research was partially supported by the National Science Foundation (DMS-1613170) and the Susan M.~Smith Professorship.} 
\and
Asaf Cohen\thanks{Department of Statistics,
	University of Haifa,
	Haifa, 31905, Israel,
	Email:
	shloshim@gmail.com,
	web: https://sites.google.com/site/asafcohentau/
}}

\date{\today}

\begin{document}

\maketitle

\begin{abstract}
\skp
We consider an $n$-player symmetric stochastic game with weak interactions between the players. Time is continuous and the horizon and the number of states are finite. We show that the value function of each of the players can be approximated by the solution of a partial differential equation called the master equation. Moreover, we analyze the fluctuations of the empirical measure of the states of the players in the game and show that it is governed by a solution to a stochastic differential equation. 
Finally, we prove the regularity of the master equation, which is required for the above results.

\skp
\noi{\bf Keywords:}    Mean-field games, master equation, fluctuations, 
  finite state control problem, Markov chains\end{abstract}

\noi{\bf AMS Classification:} 
91A06, 
60J27, 
35M99, 
93E20, 
60F05, 
91A15 
\section{Introduction}\label{sec1}
\beginsec
We consider a finite horizon game in continuous-time with weakly interacting $n$ players over a finite state space. Each player controls its own rate of transition from its state to another one, aiming to minimize a cost. The value function for each player is defined via a symmetric Nash-equilibrium. We associate the game with a partial differential equation, known as the master equations. Then we show that the solution of this equation gives an ${\cal O}(n^{-1})$ approximation to the value function. We provide a stochastic differential equation that governs the limit behavior of the fluctuations of the empirical measure of the game. 
At last, we prove the regularity of the master equation.

The theory of mean field games (MFGs) was initiated a decade ago with the independent seminal works of Lasry and Lions \cite{Lasry2006,Lasry2006b,Lasry2007}, and Huang, Malham{\'e}, and Caines \cite{Huang2006,Huang2007}. This field studies limiting models for weakly interacting $n$-player stochastic games. 
The research in this area involves the analysis of the limiting model as well as the convergence of the $n$-player games to it. The convergence was analyzed from several points of views. First, by showing that an optimal control from the MFG can be used in the $n$-player game in order to generate an asymptotic Nash equilibrium, see e.g., \cite{Cardaliaguet2013} and \cite{Carmona2013}. Lacker \cite{Lacker2015general, Lacker2017} and Fischer \cite{Fischer2017} proved a converse result: they showed that any solution of the limiting problem can be approximated by an $\eps_n$-Nash equilibrium in the prelimit game. In the influential work of Cardaliaguet, Delarue, Lasry, and Lions \cite{CDLL2015} the convergence of the value functions and empirical measure of the $n$ players to the value function and the distribution of the state of the MFG is established using the master equation. This is a parabolic partial differential equation with a terminal condition, whose variables are time, state, and measure. Its solution approximates the value function of an arbitrary player from the $n$-player game at a given time when one takes the arguments as the players' state and the empirical distribution of the other players. The master equation for finite state mean field games has been studied by Gomes, Velho, and Wolfarm in \cite{gomes2014dual, gomes2014socio}. 
A comparison between this paper and the present one is given towards the end of the introduction. For further reading on the master equation see \cite{Carmona2014} and \cite{Bensoussan2015Master}.

Continuous-time finite state MFGs were studied first by Gomes , Mohr, and Souza \cite{Gomes2013}. They showed that the value functions of the $n$- players solve a coupled system of $n$ differential equations. Recalling the symmetry between the players, it is shown further that, over a small time horizon, this system can be approximated by the solution of a coupled system of two differential equations with initial-terminal conditions, which emerges from the MFG; one of the equations stands for the value function and the other for the flow of measures. Carmona and Wang studied in \cite{Carmona2016} a version of the game that includes a major player. They find an asymptotic equilibrium in the $n$-player game using the analysis of the MFG. The authors also present the master equation in this framework, but do not use it for the approximation. More recently, Cecchin and Fischer \cite{Cecchin2017} used probabilistic methods in order to show that MFG solutions provide symmetric $\eps_n$-Nash equilibria for the $n$-player game. In \cite{BBC2017a} a finite state continuous-time $n$-player game is studied. However, the analysis is done under heavy-traffic so that the number of states increases with $n$ and the limiting problem has a diffusion noise. A numerical scheme associated with this model is studied in \cite{BBC2017b}. Another paper that studies finite state MFGs in continuous time is given in \cite{Doncel2016}. However, while the limiting problem is set in continuous-time, the prelimit problem is a discrete-time $n$-player game with a finite set of actions. As opposed to the papers mentioned in this  paragraph, after an inspiring discussion with Daniel Lacker, we perform our analysis using the master equation. 

Our game borrows the framework introduced by Gomes et al.~\cite{Gomes2013}. We provide three main results. The first one, namely Theorem \ref{thm_21}, deals with the convergence of the value functions of the $n$ players to the solution of the master equation. Specifically, we show that the average of all the individual value functions of the players 
is approximately the solution of the master equation with the same entries. The second main result, Theorem \ref{thm_22}, provides the fluctuations of the empirical measure. 
In Theorem \ref{thm1300}, adapting the results of \cite{CDLL2015} to discrete state space set-up, we show that the master equation is regular under reasonable assumptions on the data of the problem, which are required for the above results to hold. 

Let us describe some of the technical details. Recall that in general control problems, the optimal control can be expressed as a minimizer of a Hamiltonian and as so depends on the value function. Therefore, in order to establish the fluctuations, the first step is to show that the empirical measure driven by the minimizer of the Hamiltonian calculated using the value function of the representative player in the symmetric $n$-player game and the empirical measure driven by the minimizer of the Hamiltonian calculated using the master equation are close to each other up to order ${\cal O}(n^{-1})$. Then the second step is to characterize the fluctuations of the dynamics of the latter empirical distribution. To tackle the first step we couple two games. In the first game, the equilibrium strategies are driven by the value function of the representative player of the symmetric $n$-player game; and in the second game strategies are driven by the solution of the master equation. While in controlled diffusion games there is a trivial coupling: every two player with the same index (one from each game) share the same noise. A more sophisticated coupling is required in the finite state setup. 

A rigorous analysis of the fluctuations in many player diffusion games as well as the derivation of large deviation and concentration inequalities are studied  by Delarue, Lacker, and Ramanan \cite{DLR2017}. For the second step we associate every state with arrival and departure processes. By scaling them and using tightness arguments and the martingale central limit theorem and benefiting from the fact that the number of states is finite, the fluctuations can be described as a solution to a stochastic differential equation, which arguably gives more intuition about the behavior of the system. Compare that to the diffusion case, where the fluctuation limit is a solution to a stochastic partial differential equation. 
Finally, we show that a variant of the assumptions given in \cite{Gomes2013} is sufficient to get our results, which holds on an arbitrary time horizon. Our results therefore extend results from \cite{Gomes2013} to arbitrary time horizon. The results above require regularity of the master equation, which we prove in the last section by constructing it from the coupled system of forward-backward partial differential equations that were described in \cite{Gomes2013}.\footnote{During the review of this paper it was brought to our attention that Cecchin and Pelino \cite{cec-pel2017} also independently analyzed the finite-state MFG using the master equation approach. The $n$-player game in that paper is formulated using the formulation in \cite{Cecchin2017} while we follow the formulation given in \cite{Gomes2013}. We prove our main result, the fluctuations, using 
a probabilistic approach which relies on coupling, whereas \cite{cec-pel2017} uses an analytical approach relying on the convergence of the generators.  It also happens that our assumptions for the convergence results (in Section~\ref{sec2}) are slightly weaker. In obtaining sufficient conditions, Section \ref{1300} of this paper, both papers adapt the approach of \cite{CDLL2015} to discrete state space. 
}


The rest of the paper is organized as follows. In Section \ref{sec2}, we analyze the $n$-player game: we introduce the game, the master equation associated with it, and provide the convergence results. In Section \ref{sec1300}, we provide simple sufficient conditions for the results in Section \ref{sec2} to hold.

In the rest of this section we will list some frequently used notation.

\subsection{Notation}\label{sec11}
We use the following notation. For every $a,b\in\R$, $a\wedge b=\min\{a,b\}$ and $a\vee b=\max\{a,b\}$. We denote $[0,\iy)$ by $\R+$. The set of positive integers is denoted by $\N$. For every $d\in\N$ and $a,b\in\R^d$, $a\cdot b$ denotes the usual
scalar product, $\|\cdot\|$ denotes Euclidean norm, and $a^\top$ is the transpose of $a$. $\{e_1, \ldots, e_d\}$ is the standard basis of $\R^d$. Also, set $e_{xz}=e_z-e_x$. For any $k\in\N$, the set $\{1,\ldots,k\}$ is denoted by $[k]$; throughout the paper we will take $k=d,n$. 
%
For every $x\in[d]$, let $$\R^d_x:=\Big\{\eta\in\R^d:\forall [d]\ni y\ne x,
\;\eta_y\ge 0\;\;\text{and}\;\;\eta_x=-\sum_{z\in[d],z\ne x}\eta_z\Big\}$$ and  $\R^d_{[d]}=\cup_{x\in[d]}\R^d_x$. 
Let $\calQ^{d\times d}$ be the set of all $d\times d$ transition rate matrices; that is $\al\in \calQ^{d\times d}$ if for every different $x,y\in[d]$, $\al_{xy}\ge 0$, and $\sum_{z\in[d]}\al_{xz}=0$.
Let 
$\calD([0,T],\R^d)$ 
be the space of functions
that are 
right-continuous with finite left limits 
mapping $[0,T]\to\R^d$. Endow these spaces with the uniform norm topology.

\section{The $n$-player game and the MFG}\label{sec2}
\beginsec

\subsection{The stochastic model}\label{sec21} 
We consider a finite state continuous-time symmetric $n$-player game with weak interactions. The players' positions can be any of the states in $[d]$, where $d$ is an integer greater than one. At every time instant each of the players chooses the rate of transition from its own state to the others, aiming to minimize a cost. Both the transition rates and the cost depend on the current time, the player's state, and the empirical distribution of the players over the states. 

We now give a precise description of the controlled stochastic processes of interest. Fix a finite horizon $T>0$. Let $(\Omega,\calF,\{\calF_t\},\PP)$ be a filtered probability space that supports unit rate independent Poisson
processes $P^{n,i,x}$, $n\in\N, i\in[n], x\in[d]$. The controlled jump process $(X^{n,i}(t))_{t\in[0,T]}$, which stands for players $i$'s state, is taking values in $[d]$ and is defined through $P^{n,i,x}$, $x\in[d]$, in a way that $X^{n,i}$ jumps to state $x$ whenever the Cox process $\hat P^{n,i,x}$, which is introduced below, jumps. The initial states $\{X^{n,i}(0)\}_{i=1}^n$ are independent and identically distributed according to $\mu^n_0$. We assume that $\limn \mu^n_0$ exists in distribution.
Player $i$ chooses a measurable function $a^{n,i}:[0,T]\times [d]\to\R^d_{[d]}$ satisfying these following conditions: for every distinct $x,y\in[d]$, $a^{n,i}_y(t,x)\in\R_+$ and  $a^{n,i}_x(t,x)=-\sum_{y\in[d],y\ne x}a^{n,i}_y(t,x)$. For every $x\ne X^{n,i}(t)$, the expression $a^{n,i}_{x}(t, X^{n,i}(t))$ represents player $i$'s rate of transition from its current state to state $x$. 

For every $n\in\N,i\in[n]$, and $x\in[d]$, let us denote
$$\hat P^{n,i,x}(t):=P^{n,i,x}\Big(\int_0^t(a^{n,i}_x(s,X^{n,i}(s))\vee 0)ds\Big),\qquad t\in[0,T]$$
and
\begin{align}\label{201}
\rho^{n,i,x}(t):=\hat P^{n,i,x}(t)-\int_0^t(a^{n,i}_x(s,X^{n,i}(s))\vee 0)ds, \qquad t\in[0,T].
\end{align}The processes described in \eqref{201} are $\{\calF_t\}$-martingales with the predictable quadratic variations processes given by
\begin{equation}\label{202}
\langle \rho^{n,i,x},\rho^{n,j,y} \rangle(t) =\delta_{ij,xy}\int_0^t(a^{n,i}_x(s,X^{n,i}(s))\vee 0)ds,
\end{equation}
where $\delta_{ij,xy}=1$ if $i=j$ and $x=y$ and $0$ otherwise. Note that according to \eqref{201}, the Poisson process $P^{n,i,x}$ is ``not active" while $a^{n,i}_x\le0$. By our definition of the control, the last inequality obviously holds during the time the process $X^{n,i}$ equals $x$. 

The empirical distribution of the states of the players in the set $[n]\setminus\{i\}$ at time $t$ is given by 
\begin{align}\notag
\mu^{n,i}(t):=\frac{1}{n-1}\sum_{j,j\ne i}e_{X^{n,j}(t)}=\Big(\frac{1}{n-1}\sum_{j,j\ne i}\1_{\{X^{n,j}(t)=x\}}:x\in[d]\Big),
\end{align}
where hereafter the parameters $i,j$, and $k$ represent a player and belong to $[n]$. Moreover, $\sum_{j,j\ne i}$ stands for $\sum_{j\in[n],j\ne i}$.   
Also, denote
 \begin{align}\notag
\calP([d]):=\Big\{\eta\in\R^d_+ : \sum_{i\in[n]}\eta_i=1\Big\}
\end{align}
and
\begin{align}\notag
\calP^n([d]):=\{\eta\in\calP([d]): n\eta\in(\N\cup\{0\})^d\}.
\end{align}
The total expected cost for player $i$ starting at time $t$, associated with the initial condition $\\$$(X^{n,i}(t),\mu^{n,i}(t))=(x,\eta)\in[d]\times\calP^{n-1}([d])$ and the profile of strategies $a=(a^{n,1},\ldots,a^{n,n})$ is given by,
\begin{align}\notag
&J^{n,i}(t,x, \eta,a) := 
\E\Big[\int_t^T f(X^{n,i}(s),\mu^{n,i}(s), a^{n,i}(s,X^{n}(s)))ds + g(X^{n,i}(T),\mu^{n,i}(T))\Big],
\end{align}
where $f:[d]\times \calP([d])\times\R^d_{[d]}\to\R$ is the running cost and $g: [d]\times\calP([d])$ is the terminal cost. For every $(x,\eta)$, the function $a\mapsto f(x,\eta,a)$ is assumed to be independent of the $x$-th component of $a$. 
Also, $f$ and $g$ are measurable functions that will be required to satisfy additional conditions, which will be introduced in the sequel (see Assumptions \ref{assumption21} and \ref{assumption22}, and Remark \ref{rem22}). 

We are interested in finding an asymptotic 
Nash equilibrium, where the goal of each player is to minimize its own cost $J^{n,i}$ given the other player actions.
\subsection{The HJB equation}\label{sec22}

Consider the Hamiltonian 
\begin{align}\notag
H(x,\eta,p)=\inf_{a\in\R^d_x}h(x,\eta,a,p),
\end{align}
where $h:[d]\times\calP([d])\times\R^d_x\times\R^d\to\R$ is defined as
\begin{align}\notag
h(x,\eta,a,p)=f(x,\eta,a)+\sum_{y,y\ne x}a_yp_y.
\end{align}
Throughout, the parameters $x,y$, and $z$ represent a state of the system and belong to $[d]$. Also, $\sum_{y,y\ne x}$ stands for $\sum_{y\in[d], y\ne x}$.  
The following assumption requires a unique minimizer for the Hamiltonian. 
\begin{assumption}\label{assumption21}$\,$
There is a unique measurable function 
$ a^*:[d]\times \calP([d])\times \R^d\to\R_{[d]}^d$ such that
\begin{align}\notag
a^*(x,\eta,p)=\underset{a\in \R_x^d}{\arg\min}\;h(x,\eta,a,p),
\end{align}
where uniqueness is considered up to its $x$-th coordinate, since $f(x,\eta,a)$ is independent of the $x$-th component of $a$.
\end{assumption}
We set $a^*_x(x,\eta,p)=-\sum_{y,y\ne x }a^*_y(x,\eta,p)$. This assumption enables us to establish the uniqueness of the optimal response of any reference player to the other players' strategies, and in the MFG, which will be introduced in the next subsection, the uniqueness of the optimal control of the decision maker. It is a standard assumption in the MFG literature, see \cite[(5)]{Gomes2013}, \cite[Hypotheses 8]{Carmona2016}, and \cite[Assumption (C)]{Cecchin2017} in the context of finite state MFGs and \cite[Lemma 2.1]{Carmona2013} \cite[(13)]{CDLL2015}, and \cite[Assumption 3.1.(b)]{BBC2017a} for diffusion scaled MFGs.

The differential equations associated with this game are given in \cite{Gomes2013} and in \cite{Carmona2016}. However, for completeness of the presentation, we shortly illustrate its derivation. Fix $i\in[n]$. We refer to player $i$ as a {\it reference} player. Assume that all the players in the set $[n]\setminus\{i\}$ 
use the same control $a^{n,-i}$. Denote the cost for player $i$ who uses the control $a$ in this case by $J^{n,i}_{a^{n,-i}}(t,x, \eta,a)$. Also, set
\begin{align}\notag
V^{n,i}(t,x,\eta)=V^{n,i}_{a^{n,-i}}(t,x,\eta):=\inf_{a^{n,i}} J^{n,i}_{a^{n,-i}}(t,x,\eta,a^{n,i}),\;(t,x,\eta)\in[0,T]\times[d]\times\calP^{n-1}([d]),
\end{align}
where the infimum is taken over controls $a^{n,i}$.
Then $V^{n,i}$ solves
\begin{align}\notag
\begin{cases}
-\partial_t V^{n,i}(t,x,\eta)=\underset{y,z}{\sum}\Big[V^{n,i}\left(t,x,\eta+\tfrac{1}{n-1}e_{yz}\right)-V^{n,i}(t,x,\eta)\Big](n-1)\eta_ya^{n,-i}_z(t,y)\\
\qquad\qquad\qquad\qquad+H\left(x,\eta,\Delta_xV^{n,i}(t,\cdot,\eta)\right),\\
V^{n,i}(T,x,\eta)=g(x,\eta),
\end{cases}
\end{align}
where $(t,x,\eta)\in[0,T]\times[d]\times\calP^{n-1}([d])$ and for every $x\in[d]$ and $\phi:[d]\to\R$, 
\begin{align}\notag
\Delta_x\phi(\cdot):=\left(\phi(l)-\phi(x)
: l\in[d]\right),
\end{align}
and $\sum_{y,z}=\sum_{y\in[d]}\sum_{z\in[d]}$.

In a symmetric equilibrium we expect that the controls of all the players would be the same as the one of player $i$. That is, each player $j\in[n]$ would use the same control:
$$
a^*_z(X^{n,j}(t),\mu^{n,j}(t),\Delta_{X^{n,j}(t)}V^{n,j}(t,\cdot,\eta)),\qquad (t,z)\in[0,T]\times[d].
$$
The symmetry between the players also yields that $V^{n,j}=V^{n,i}$ for every $j\in[n]$. We denote $V^n=V^{n,i}$. 
Therefore, the last system becomes
\begin{align}\tag{{\bf \text{HJB}}\bf{(n)}}\label{212}\\\notag
\begin{cases}
-\partial_t V^{n}(t,x,\eta)=\underset{y,z}{\sum}\eta_yD^{n,y,z}V^{n}(t,x,\eta)a^*_z(y,\eta+\tfrac{1}{n-1}e_{yx},\Delta_yV^{n}(t,\cdot,\eta+\tfrac{1}{n-1}e_{yx}))\\
\qquad\qquad\qquad\quad\;\;+H\left(x,\eta,\Delta_xV^{n}(t,\cdot,\eta)\right),\\
V^{n}(T,x,\eta)=g(x,\eta),
\end{cases}
\end{align}
where\footnote{For simplicity, we are abusing notation here and in the sequel and write $\sum_{y,z}$ for $\sum_{y\in[d],\eta(y)\ne 0}\sum_{z\in[d]}$ whenever the coefficient $\eta_y$ appears in the sum.}
\begin{align}\label{213}
D^{n,y,z}\phi(t,x,\eta):=(n-1)\left(\phi\left(t,x,\eta+\tfrac{1}{n-1}e_{yz}\right)-\phi(t,x,\eta)\right),
\end{align}
for $(t,x,\eta)\in[0,T]\times[d]\times\calP^{n-1}([d])$. 
The reason for the additional term $\tfrac{1}{n-1}e_{yx}$ in the argument of $a^*$ 
is that we need to modify $\eta$ so that it represents the position of the player that stands for position $y$. 

Gomes et al.~(2013) provide sufficient conditions for Assumption \ref{assumption21} to hold and for the existence and uniqueness of a solution of this system, in addition to equilibrium strategies, which we now summarize.

\begin{lemma}\label{lem21}[Proposition 1 and Theorems 4, 5, and 6 in \cite{Gomes2013}]
Suppose that the running cost $f$ and the terminal cost $g$ satisfy the following conditions.
\begin{enumerate}
\item (Lipschitz-continuity) The function $f$ is differentiable with sepect to (w.r.t.) $a$ and there exists a constant $c_L>0$ such that for every $x\in[d]$, $\eta,\eta'\in\calP([d])$, and $a\in\R^d_+$,
\begin{align}\notag
&|g(x,\eta)-g(x,\eta')|+|f(x,\eta,a)-f(x,\eta',a)|+\|\nabla_af(x,\eta,a)-\nabla_af(x,\eta',a)\|\\\notag
&\qquad \le c_L\|\eta-\eta'\|.
\end{align}
\item (strong convexity) There exists a constant $c_V>0$ such that for every $x\in[d]$, $\eta\in\calP([d])$, and $a,a'\in\R^d_+$,
\begin{align}\notag
f(x,\eta,a)-f(x,\eta,a')\ge \nabla_af(x,\eta,a')(a-a') + c_V\|a-a'\|^2.
\end{align}
\item (superlinearity) For every $(x,\eta)\in[d]\times\calP([d])$ and $\{a_y\}_{y\ne x}\subset\R_+$, 
\begin{align}\label{215-a}
\lim_{a_x\to\iy}\frac{f(x,\eta,a)}{\|a\|}=\iy.
\end{align}
\end{enumerate}
Then, 
\begin{enumerate}
\item[a.] Assumption \ref{assumption21} holds and the function $a^*$ is uniformly Lipschitz in its arguments. 
\item[b.] The Hamiltonian $H$ is uniformly Lipschitz in its arguments over compact sets.
\item[c.] 
The system \eqref{212} has a unique classical solution, which coincides with the value function associated with equilibrium strategies that for every $j\in\N$ are given by $$a^{n,j}(t,X^{n,j}(t))=a^*(X^{n,j}(t),\mu^{n,j}(t),\Delta_{X^{n,j}(t)}V^{n}(t,\cdot,\mu^{n,j}(t))).$$
\end{enumerate}
\end{lemma}

\subsection{The master equation and the MFG}\label{sec23}
In this section we present the master equation associated with the game described in the previous subsection. We provide heuristics for its derivation from both PDE and probabilistic perspectives. The MFG is fully analyzed in \cite{Gomes2013}. A rigorous derivation of the master equation from the MFG is given in Section \ref{sec1300}. In the next subsection we focus on the relationship between the master equation and \eqref{212}, which is used to establish the fluctuations of the empirical measure in Section \ref{sec25}. 

The master equation associated with the $n$-player game is given by
\begin{align}\tag{{\bf ME}}\label{215}
\begin{cases}
-\partial_t U(t,x,\eta)=\sum_{y,z}\partial_{\eta_z}U(t,x,\eta)\eta_ya^*_z(y,\eta,\Delta_yU(t,\cdot,\eta))+H(x,\eta,\Delta_{x}U(t,\cdot,\eta)),\\
U(T,x,\eta)=g(x,\eta),
\end{cases}
\end{align}
where $U:[0,T]\times [d]\times\calP([d])\to\R$.
Informally, the structure of the master equation follows from \eqref{212} using the approximation 
\begin{align}\label{216}
&\left(V^{n}\left(t,x,\eta+\tfrac{1}{n-1}e_{yz}\right)-V^{n}(t,x,\eta)\right)(n-1)\approx      \partial_{\eta_z} V^{n}(t,x,\eta)- \partial_{\eta_y} V^{n}(t,x,\eta),
 \end{align} 
together with our setting $\sum_{z\in[d] }a^*_z(y,\eta,p)=0$, $(z,\eta,p)\in[d]\times\calP([d])\times\R^d$.

We now provide a probabilistic illustration of the master equation using a MFG problem. The MFG was studied by Gomes et al.~in \cite{Gomes2013} and \cite{Cecchin2017} and with an additional major player in \cite{Carmona2016}. 
Fix $\eta\in\calP([d])$ and a function $\beta:[0,T]\to\calQ^{d\times d}$, see Section \ref{sec11}. 
Using the terminologies from the $n$-player game, we consider a reference player, which we refer to as the {\it decision maker} (DM). The matrix $\beta$ stands for the transition matrix of the non-reference players. Let $\mu^\beta:[0,T]\to\calP([d])$ be the solution of 
\begin{align}\label{217}
\begin{cases}
\frac{d}{dt}\mu^\beta_y(t)=\sum_{x\in[d]}\mu^\beta_x(t)\beta_{xy}(t),\quad (t,y)\in[0,T]\times[d],\\
\mu^\beta(0)=\mu_0.
\end{cases}
\end{align}
The function $\mu^\beta$ represents the limit of the evolution in time of the distribution of the non-reference players. Consider the following Markov decision problem. At any time instant, the DM chooses the rate of transition from its current state to another state in the set $[d]$ as follows. Let $(X(t))_{t\in[0,T]}$ be a $[d]$-valued controlled process that stands for the state process of the DM. A Markovian control of the DM is a process $(a(t))_{t\in[0,T]}$ taking values in $\R^d_{[d]}$, such that for every $t\in[0,T]$ and every $x\ne X(t)$, $a_x(t)$ is the transition rate from the current state $X(t)$ to state $x$. Also, set $a_{X(t)}(t)=-\sum_{x\ne X(t)}a_x(t)$. The total
expected cost for the DM starting at time $t$, associated with the initial condition $(X(t),\mu^\beta(t))=(x,\eta)\in[d]\times\calP([d])$ and the control $a$ is given by,
\begin{align}\notag
&J_\beta(t,x, \eta,a) := 
\E\Big[\int_t^T f(X(s),\mu^\beta(s),a(s))ds + g(X(T),\mu^\beta(T))\Big].
\end{align}
The value function is therefore given by
\begin{align}\notag
V_\beta(t,x,\eta)=\inf_{a} J_\beta(t,x,\eta,a),\qquad(t,x,\eta)\in[0,T]\times[d]\times\calP([d]),
\end{align}
where the infimum is taken over a set of admissible controls, which we do not define rigorously. The value function solves the following HJB
\begin{align}\label{220}
\begin{cases}
-\frac{d}{dt}V_\beta(t,x,\mu^\beta(t))=H(x,\mu^\beta(t),\Delta_xV_\beta(t,\cdot,\mu^\beta(t))),\qquad (t,k)\in[0,T]\times[d],\\
V_\beta(T,x,\mu^\beta(T))=g(x,\mu^\beta(T)),
\end{cases}
\end{align}
At the end of this subsection we provide sufficient conditions for the existence and uniqueness of a solution to the system given in \eqref{217} and \eqref{220}. 
This system gives,
\begin{align}\notag
-\partial_tV_\beta(t,x,\mu^\beta(t))&=\sum_{z\in[d]}\partial_{\eta_z}V_\beta(t,x,\mu^\beta(t))\frac{d}{dt}\mu^\beta_z(t)+H(x,\Delta_xV_\beta(t,\cdot,\mu^\beta(t)))\\\notag
&=\sum_{y,z}\partial_{\eta_z}V_\beta(t,x,\mu^\beta(t))\mu^\beta_z(t)\beta_{yz}(t)+H(x,\Delta_xV_\beta(t,\cdot,\mu^\beta(t))).
\end{align}
The optimal control is thus given by
\begin{align}\notag
a^\beta(t)=a^*(X(t),\mu^\beta(t),\Delta_{X(t)}V_\beta(t,\cdot,\mu^\beta(t))),\qquad t\in[0,T].
\end{align}
In a symmetric equilibrium all players are expected to use $\beta$ that satisfies
$$
\beta_{xy}(t)=a^*_y(x,\mu^\beta(t),\Delta_xV_\beta(t,\cdot,\mu^\beta(t))),\qquad(t,x,y)\in[0,T]\times[d]^2.
$$
Hence, the forward-backward system is

\begin{align}
\begin{cases}\label{224united}
\frac{d}{dt}\mu_y(t)=\sum_{x\in[d]}\mu_x(t)a^*_y(x,\mu(t),\Delta_xV(t,\cdot,\mu(t))),\\
-\partial_tV(t,x,\mu(t))=\sum_{y,z}\partial_{\eta_z}V(t,x,\mu(t))\mu_z(t)a^*_z(x,\mu(t),\Delta_xV(t,\cdot,\mu(t)))\\
\qquad\qquad\qquad\quad+H(x,\mu(t),\Delta_xV(t,\cdot,\mu(t))),\\
\mu(0)=\mu_0, \qquad V(T,x)=g(x,\mu(T)),
\end{cases}
\end{align}
where $(t,x,y,\mu_0)\in[0,T]\times[d]^2\times\calP([d])$. Observe that the master equation comes from the second differential equation. We will make this connection more rigorous in Section~\ref{sec1300}, where we construct a smooth solution of the master equation using the forward-backward system.
%

The next lemma tells us when the above system has a unique smooth solution.
\begin{lemma}\label{lem22}[Proposition 4 and Theorem 2 in \cite{Gomes2013}]
Under the assumptions from Lemma \ref{lem21} there exists a solution to the system of equations given in \eqref{224united}. Moreover, under the following additional assumptions we also have uniqueness.
\begin{enumerate}
\item For any $\eta,\eta'\in\calP([d])$, 
\begin{align}\notag
\sum_{x\in[d]}(\eta_x-\eta_x')(g(x,\eta)-g(x,\eta'))\ge 0.
\end{align}
\item
For every $M>0$ there exists a parameter $c_M>0$ such that for every $x\in[d]$, $\eta\in\calP([d])$, and $p,p'\in[-M,M]^d$,
 \begin{align}\label{224y}
H(x,\eta, p)-H(x,\eta, p')-a^*(x,\eta,p')\cdot(p-p')\le -c_M\|p-p'\|^2.
\end{align}
\item There exists a positive parameter $c_H$ such that for every $t\in[0,T]$, $\eta,\eta'\in\calP([d])$, and $p,p'\in\R^d$,
 \begin{align}\notag
\eta\cdot(H(\cdot,\eta',p)-H(\cdot,\eta,p))+\eta'\cdot(H(\cdot,\eta,p')-H(\cdot,\eta',p'))\le-c_H\|\eta-\eta'\|^2,
\end{align}
where $H(\cdot,\eta,p):=(H(x,\eta, p):x\in[d])$. 
\end{enumerate}
\end{lemma}

\begin{remark}\label{rem21} [Section 2.8 in \cite{Gomes2013}]
Simple sufficient conditions for the assumptions above are that $\eta\mapsto (g(x,\eta):x\in[d])$ is the gradient of a convex function and that $H$ can be rewritten as 
\begin{align}\notag
H(x,\eta,\Delta_xz)=H_1(x,\Delta_xz)+H_2(x,\eta),
\end{align}
with $H_1$ satisfying \eqref{224y} and $H_2$ satisfying 
\begin{align}\notag
(H_2(\cdot, \eta')-H_2(\cdot,\eta))\cdot(\eta'-\eta)\ge c_H\|\eta-\eta'\|^2,
\end{align}
where $H_2(\cdot,\eta):=(H_2(x,\eta):x\in[d])$. 
The last property holds, for instance, if $\eta\mapsto (H_2(x,\eta):x\in[d])$ is the gradient of a convex function.
\end{remark}

\subsection{Convergence of the value functions} \label{sec24}
We start by introducing key conditions that are necessary for all the results in the rest of Section \ref{sec2}. As we state below, some of the conditions are derived from the conditions given in Lemma \ref{lem21}. Yet, we choose to state only the most primitive conditions that are necessary for the proofs. 
\begin{assumption}\label{assumption22}
\begin{enumerate}
\item The systems of the differential equations \eqref{212} and \eqref{215} have unique classical solutions.
\item 
The function $\nabla_\eta U$ is uniformly bounded over $[0,T]\times[d]\times\calP([d])$. Moreover, 
There exists $c_L>0$ such that for every $t\in[0,T]$, $x\in[d]$, and $\eta,\eta'\in\calP([d])$,
\begin{align}\label{225-c}
\|\nabla_\eta U(t,x,\eta)-\nabla_\eta U(t,x,\eta')\|\le c_L\|\eta-\eta'\|.
\end{align}
\end{enumerate}
\noi Denote $c_U:=\sup_{(t,x,\eta)}|U(t,x,\eta)|$, which is finite by the continuity of $U$ over its compact domain.
\begin{enumerate}[resume]
\item Modifying $c_L$ if necessary, we assume\footnote{In fact for the result in this subsection we do not need the Lipschitz-continuity of the control w.r.t.~$\eta$. It is used to show the fluctuations in Section \ref{sec25}. We choose to state it here for the presentation's sake.} that for every $M>0$ there is a parameter $c_M>0$ such that for every $x\in[d]$, $\eta,\eta'\in\calP([d])$, and $p,p'\in[-M,M]^d$,
\begin{align}\label{225-a}
\|a^*(x,\eta,p)-a^*(x,\eta',p')\|\le c_M\left(\|\eta-\eta'\|+\|p-p'\|\right).
\end{align}
Moreover, for every distinct $x,y\in[d]$, the function $a^*_y(x,\cdot,\cdot)$ is uniformly bounded below away from zero.

\item For every $M>0$ there is a parameter $c_M>0$ such that for every $x\in[d]$, and $\eta\in\calP([d])$, and $p,p'\in[-M,M]^d$,
\begin{align}\notag
|H(x,\eta,p)-H(x,\eta,p')|\le c_M\|p-p'\|.
\end{align}
\end{enumerate}
\end{assumption}
\begin{remark}\label{rem22}
i. 
Recall that the existence and uniqueness of a solution to \eqref{212} is guaranteed under the conditions from Lemma \ref{lem21}. The existence and uniqueness of a solution to \eqref{215} as well as its properties are given in Section \ref{sec1300}. From \cite[Proposition 6]{Gomes2013} the value functions $V^{n,i}$ are bounded above, uniformly in $(t,x,\eta,n,i)$. Set $c_V:=c_U\vee\sup_{(t,x,\eta,n,i)}|V^{n,i}(t,x,\eta)|$, Notice that the continuity of $a^*$ over the compact domain $[d]\times\calP([d])\times[-2c_V,2c_V]$ implies that it is uniformly bounded above. Together with the third part of the assumption we get that there is a constant $c_B>0$ such that for every $x,y\in[d], \eta\in\calP([d])$, and $p\in[-2c_V,2c_V]^d$, 
\begin{align}\label{COHEN001}
\frac{1}{c_B}\le a^*_y(x,\eta,p)\le c_B.
\end{align}
The uniform lower bound of the optimal appears in \cite[Definition 3.2.(iii)]{DRW2016} and \cite[Lemma 3]{Cecchin2017}. It is somewhat equivalent to the requirement in the diffusion case that the diffusion coefficient, which drives the noise, is bounded away from zero. Such a condition holds if we restrict the controls to be greater than a given positive constant, or alternatively, if the running cost satisfies \eqref{215-a} with $a_x\to 0^+$.

ii. 
In the light of \eqref{216}, the uniform Lipschitz-continuity of 
$\nabla_\eta U$ is needed below in order to uniformly approximate $ D^{n,y,z}U$ by $\partial_{\eta_z}U$. 
Kolokoltsov and Yang (2013) provide sufficient conditions for a statement equivalent to Assumption \ref{assumption22}.2 in the diffusion case, see \cite[Theorem 4.3]{Kolokoltsov2013}. The Lipschitz-continuity of $a^*$ and $H$ is a standard requirement in the literature of MFGs. As stated in Lemma \ref{lem21} above, Gomes et al.~(2013) derives them as well as Assumption \ref{assumption21} from their assumptions. For further references in the context of finite state MFGs see \cite[Hypothesis 3]{Carmona2016} and \cite[Section 2.1.1]{Cecchin2017}. For the references in the diffusion setup see e.g., \cite[Assumptions (A.2), (A.3), (A.4), and (A.5)]{Carmona2013}, \cite[Assumption (A3)]{Kolokoltsov2013}, \cite[Section 2.3]{CDLL2015}, and \cite[Assumptions 3.1.(a) and 5.1]{BBC2017a}. 

iii. The convergence results of Gomes et al.~\cite{Gomes2013} hold only if the horizon $T$ is sufficiently small. This restriction emerges from their need that the functions $\{V^n\}_n$ be uniformly Lipschitz, see \cite[Proposition 7]{Gomes2013}. Using the master equation to approximate the value function, we bypass this requirement and therefore can handle an arbitrary time horizon. Specifically, instead of working with the equilibrium controls that can be expressed using $V^n$, we use controls that are defined using the solution of the master equation $U$.  This is because the empirical distribution generated by the equilibrium controls can be approximated by the empirical distribution of the states generated by the controls which are defined through $U$, as we demonstrate.
As a result, in order to attain the convergence of $V^n$ to $U$ we only need that $\nabla_\eta U$ is Lipschitz as is shown in Section \ref{sec1300}. 
\end{remark}

    
Consider the functions $V^{n}$ and $U$ that solve \eqref{212} and \eqref{215}, respectively.
The next theorem provides an ${\cal O}(n^{-1})$ approximation to the value function in the $n$-player game. 
An equivalent result in the diffusion case was established in \cite[Theorem 6.6]{CDLL2015}.
\begin{theorem}\label{thm_21}
Set $W^n:=V^{n}-U$. Under Assumptions \ref{assumption21} and \ref{assumption22}, there exists $C>0$ such that for every $n\ge 2$ 
it holds that
\begin{align}\notag
&\sup_{t\in[0,T]} \E\left[\left|W^n\left(t,X^{n,i}(t),\mu^{n,i}(t)\right)\right|^2\right]\\\label{225}
&\qquad+\int_{0}^T\E\left[\left\|\Delta_{X^{n,i}(s)}W^n\left(s,X^{n,i}(s),\mu^{n,i}(s)\right)\right\|^2\right]ds\le\frac{C}{n^2}.
\end{align}
Moreover, setting
${\bf x}=(x_1,\ldots,x_n)\in[d]^n$, and defining $m_{\bf x}^{n,i}:=\frac{1}{n-1}\sum_{j,j\ne i}e_{x_j}$, we have that
\begin{align}\label{224z}
\sup_{t\in[0,T]}\frac{1}{n}\sum_{i\in[n]}\left|V^{n}(t,x_i,m_{\bf x}^{n,i})-U(t,x_i,m_{\bf x}^{n,i})\right|\le\frac{C}{n}.
\end{align}  
\end{theorem}
%
{\bf Proof.}   
Inequality \eqref{224z} follows from \eqref{225}, which we now prove. 
Fix $i\in[n]$. Let us denote for every $(t,y,j)\in[0,T]\times[d]\times[n]$,
\begin{align}\notag
a^{n,j}_y(t)=a^*_y\left(X^{n,j}(t),\mu^{n,j}(t),\Delta_{X^{n,j}(t)}V^{n}(t,\cdot,\mu^{n,j}(t)\right),
\end{align}
and for every $t\in[0,T]$ and $z\in[d]$,
\begin{align}\notag
B_{i,z}(t)
&=\left|W^n\left(t,z,\mu^{n,i}(t)\right)\right|^2-\left|W^n\left(t,X^{n,i}(t),\mu^{n,i}(t)\right)\right|^2\end{align}
and for every $k\ne i$,
\begin{align}\notag
B_{k,z}(t)&=\left|W^n\left(t,X^{n,i}(t),\mu^{n,i}(t)+\tfrac{1}{n-1}e_{X^{n,k}(t)z}\right)\right|^2-\left|W^n\left(t,X^{n,i}(t),\mu^{n,i}(t)\right)\right|^2\\\notag
&=\left[W^n\left(t,X^{n,i}(t),\mu^{n,i}(t)+\tfrac{1}{n-1}e_{X^{n,k}(t)z}\right)+W^n\left(t,X^{n,i}(t),\mu^{n,i}(t)\right)\right]\\\notag
&\quad\times \frac{1}{n-1}D^{n,X^{n,k}(t),z}W^n\left(t,X^{n,i}(t),\mu^{n,i}(t)\right).
\end{align}
Recall the definitions of $\hat P^{n,j,x}$ and $\rho^{n,j,x}$ from \eqref{201} and the paragraph preceding it. 
Now,
applying It\^{o}'s formula to $|W^n|^2$ and using \eqref{212} and \eqref{215} we get that 
\begin{align}\notag
&\left|W^n\left(T,X^{n,i}(T),\mu^{n,i}(T)\right)\right|^2-\left|W^n\left(t,X^{n,i}(t),\mu^{n,i}(t)\right)\right|^2
\\\label{228}
&=2\int_t^TW^n\left(s,X^{n,i}(s),\mu^{n,i}(s)\right)\partial_tW^n\left(s,X^{n,i}(s),\mu^{n,i}(s)\right)ds\\\notag
&\quad+\int_t^T\Big\{\sum_{j\in[n]}\sum_{z\ne X^{n,j}(s)}B_{j,z}(s-)d\hat P^{n,j,z}(s)\Big\}\\\notag
&=\Upsilon+\int_t^T\left\{
\Phi_1(s)+\Phi_2(s)
+2W^n\left(s,X^{n,i}(s),\mu^{n,i}(s)\right)[\Psi_1(s)+\Psi_2(s)+\Psi_3(s)+\Psi_4(s)]\right\}ds,
\end{align}
where
\begin{align}\notag
\Upsilon&:=\int_t^T\sum_{j\in[n]}\sum_{z\ne X^{n,j}(s)}B_{j,z}(s-)d\rho^{n,j,z}(s),\\\notag
\Phi_1(s)&:=\sum_{j,j\ne i}\sum_{z\ne X^{n,j}(s)}B_{j,z}(s)a^{n,j}_z(s),\\\notag
\Phi_2(s)&:=\sum_{z\ne X^{n,i}(s)}B_{i,z}(s)a^{n,i}_z(s),
\end{align}
and
\begin{align}
\Psi_1(s)&:=\sum_{y,z}\Big[a^*_z\left(y,\mu^{n,i}(s),\Delta_yU\left(s,\cdot,\mu^{n,i}(s)\right)\right)\\\notag
&\qquad\quad-a^*_z\left(y,\mu^{n,i}(s)+\tfrac{1}{n-1}e_{yX^{n,i}(s)},\Delta_yV^{n}\left(s,\cdot,\mu^{n,i}(s)+\tfrac{1}{n-1}e_{yX^{n,i}(s)}\right)\right)\Big]\\\notag
&\qquad\qquad\times\mu^{n,i}_y(s)\partial_{\eta_z}U\left(s,X^{n,i}(s),\mu^{n,i}(s)\right),\\\notag
\Psi_2(s)&:=\sum_{y,z}\left(\partial_{\eta_z}-D^{n,y,z}\right)U\left(s,X^{n,i}(s),\mu^{n,i}(s)\right)\mu^{n,i}_y(s)\\\notag
&\qquad\quad\times a^*_z\left(y,\mu^{n,i}(s)+\tfrac{1}{n-1}e_{yX^{n,i}(s)},\Delta_yV^{n}\left(s,\cdot,\mu^{n,i}(s)+\tfrac{1}{n-1}e_{yX^{n,i}(s)}\right)\right),\\\notag
\Psi_3(s)&:=-\sum_{y,z}D^{n,y,z}W^n\left(s,X^{n,i}(s),\mu^{n,i}(s)\right)\mu^{n,i}_y(s)\\\notag
&\qquad\quad\times a^*_z\left(y,\mu^{n,i}(s)+\tfrac{1}{n-1}e_{yX^{n,i}(s)},\Delta_yV^{n}\left(s,\cdot,\mu^{n,i}(s)+\tfrac{1}{n-1}e_{yX^{n,i}(s)}\right)\right),\\\notag
\Psi_4(s)&:=\left[H\left(X^{n,i}(s),\mu^{n,i}(s),\Delta_{X^{n,i}(s)}U\left(s,\cdot,\mu^{n,i}(s)\right)\right)\right.\\\notag
&\qquad\qquad\left.-H\left(X^{n,i}(s),\mu^{n,i}(s),\Delta_{X^{n,i}(s)}V^{n}\left(s,\cdot,\mu^{n,i}(s)\right)\right)\right].
\end{align}

Notice that 
\begin{align}\notag
&\Phi_1(s)+2W^n\left(s,X^{n,i}(s),\mu^{n,i}(s)\right)\Psi_3(s)\\\notag
&\quad=
\frac{1}{n-1}\sum_{y,z}\left|D^{n,y,z}W^n\left(s,X^{n,i}(s),\mu^{n,i}(s)\right)\right|^2\mu^{n,i}_y(s)\\\notag
&\qquad\quad\times a^*_z\left(y,\mu^{n,i}(s)+\tfrac{1}{n-1}e_{yX^{n,i}(s)},\Delta_yV^{n}\left(s,\cdot,\mu^{n,i}(s)+\tfrac{1}{n-1}e_{yX^{n,i}(s)}\right)\right)\\\notag
&\quad\ge0
\end{align}
and
\begin{align}\notag
\Phi_2(s)&=\sum_{z\ne X^{n,i}(s)}\left(W^n\left(s,X^{n,i}(s),\mu^{n,i}(s)\right)-W^n\left(s,z,\mu^{n,i}(s)\right)\right)^2a^{n,i}_z(s)\\\notag
&\qquad+2W^n\left(s,X^{n,i}(s),\mu^{n,i}(s)\right)\Psi_5(s),
\end{align}
where we set
\begin{align}\notag
\Psi_5(s):=\sum_{z\ne X^{n,i}(s)}\left[W^n\left(s,z,\mu^{n,i}(s)\right)
-
W^n\left(s,X^{n,i}(s),\mu^{n,i}(s)\right)
\right]
a^{n,i}_z(s).
\end{align}
By Cauchy--Schwartz inequality and \eqref{COHEN001},
\begin{align}\notag
|\Psi_5(s)|\le\sqrt{d}c_B\left\|\Delta_{X^{n,i}(s)}W^n\left(s,\cdot,\mu^{n,i}(s)\right)\right\|.
 \end{align} 
Another way of presenting $\Psi_1$ is given by
\begin{align}\notag
&\Psi_1(s)=\frac{1}{n}\sum_{j\in[n]}\sum_{z\in[d]}
\partial_{\eta_z}U\left(s,X^{n,i}(s),\mu^{n,i}(s)\right)\\\notag
&\qquad\qquad\qquad\qquad\times\Big[a^*_z\left(X^{n,j}(s),\mu^{n,i}(s),\Delta_{X^{n,j}(s)}U\left(s,\cdot,\mu^{n,i}(s)\right)\right)\\\notag
&\qquad\qquad\qquad\qquad\qquad-a^*_z\left(X^{n,j}(s),\mu^{n,j}(s),\Delta_{X^{n,j}(s)}V^{n}\left(s,\cdot,\mu^{n,j}(s)\right)\right)\Big].
\end{align}
Recalling Assumptions \ref{assumption22}.2 and \ref{assumption22}.3, Remark \ref{rem22}, and noticing that for sufficiently large $n$, $\\$$\sup_{s\in[0,T]}|\mu^{n,i}(s)-\mu^{n,j}(s)|\le 3n^{-1}$ we get that 
\begin{align}\notag
|\Psi_1(s)|
&\le\frac{C}{n}\sum_{j\in[n]}\left\|\Delta_{X^{n,j}(s)}W^n\left(s,\cdot,\mu^{n,j}(s)\right)\right\|+\frac{C}{n},
\end{align}
where in the above expression and in the rest of the proof, $C$ refers to a finite positive constant that is independent of $t$ and $n$ and which can change from one line to the next.
From Assumption \ref{assumption22}.3, $a^*$ is uniformly bounded. Moreover,  Assumption \ref{assumption22}.2 and Lemma \ref{lem_a1}, imply the uniform bound $|\left(\partial_{\eta_z}-D^{n,y,z}\right)U|\le C/n$ and thus,
$|\Psi_2(t)|\le \frac{C}{n}$. Finally, Assumption \ref{assumption22}.4 implies,
\begin{align}\notag
|\Psi_4(t)|&\le C\left\|\Delta_{X^{n,i}(t)}W^n\left(t,\cdot,\mu^{n,i}(t)\right)\right\|.
\end{align}

From the martingale property of $\rho^{n,i,z}$'s, $\E[\Upsilon]=0$. Taking expectation on both sides of \eqref{228} and recalling that $W^n(T,X^{n,i}(T),\mu^{n,i}(T))=0$, \eqref{COHEN001}, and the above equalities and estimates yield
\begin{align}\label{COHEN555}
& \E\left[\left|W^n\left(t,X^{n,i}(t),\mu^{n,i}(t)\right)\right|^2\right]+\frac{1}{c_B}\int_{t}^T\E\left[\left\|\Delta_{X^{n,i}(s)}W^n\left(s,\cdot,\mu^{n,i}(s)\right)\right\|^2\right]ds\\\notag
&\quad\le \E\left[\left|W^n\left(t,X^{n,i}(t),\mu^{n,i}(t)\right)\right|^2\right]\\\notag
&\qquad+\int_{t}^T\E\Big[\sum_{z\ne X^{n,i}(s)}\left(W^n\left(s,X^{n,i}(s),\mu^{n,i}(s)\right)-W^n\left(s,z,\mu^{n,i}(s)\right)\right)^2a^{n,i}_z(s)\Big]ds\\\notag
&\quad\le
2\E\Big[\int_t^T
\left|W^n\left(s,X^{n,i}(s),\mu^{n,i}(s)\right)\right|\Big(\left|\Psi_1(s)\right|+\left|\Psi_2(s)\right|+\left|\Psi_4(s)\right|+\left|\Psi_5(s)\right|\Big)ds\\\notag
&\quad\le
C_1\E\Big[\int_t^T
\left|W^n\left(s,X^{n,i}(s),\mu^{n,i}(s)\right)\right|\Big(\frac{1}{n}\sum_{j\in[n]}\left\|\Delta_{X^{n,j}(s)}W^n\left(s,\cdot,\mu^{n,j}(s)\right)\right\|
\\\notag
&\qquad\qquad\quad\qquad\qquad\qquad\qquad\qquad\qquad\qquad+
\left\|\Delta_{X^{n,i}(s)}W^n\left(s,\cdot,\mu^{n,i}(s)\right)\right\|
+
\frac{1}{n}\Big)ds,
\end{align}
where $C_1>0$ and is independent of $t$ and $n$. 
 Applying Young's inequality, $|uv|\le \eps u^2/2+v^2/(2\eps)$, for all $u,v\in\R$ and $\eps>0$, to the $n+2$ products that we obtained in the last line of the equation separately with $\eps=(2c_BC_1)^{-1}$, we get 
\begin{align}\notag
&\le 3C_1^2c_B\E\Big[\int_t^T
\left|W^n\left(s,X^{n,i}(s),\mu^{n,i}(s)\right)\right|^2ds\Big]
+
\frac{1}{4c_Bn}\int_t^T\E\Big[\sum_{j\in[n]}\left\|\Delta_{X^{n,j}(s)}W^n\left(s,\cdot,\mu^{n,j}(s)\right)\right\|^2\Big]ds
\\\notag
&\quad+\frac{1}{4c_B}\int_t^T\E\Big[\left\|\Delta_{X^{n,i}(s)}W^n\left(s,\cdot,\mu^{n,i}(s)\right)\right\|^2\Big]ds+\frac{T}{4c_Bn^2}.
\end{align}
Recall that the random variables $\{X^{n,j}(t)\}_{j=1}^n$ and therefore also the random variables $\\$$\left\{\Delta_{X^{n,j}(t)}W^n\left(t,\cdot,\mu^{n,j}(t)\right):j\in[n]\right\}$ are identically distributed. Hence, we obtain the following bound
\begin{align}\notag
&\le 3C_1^2c_B\E\Big[\int_t^T
\left|W^n\left(s,X^{n,i}(s),\mu^{n,i}(s)\right)\right|^2ds\Big]
+
\frac{1}{2c_B}\int_t^T\E\Big[\left\|\Delta_{X^{n,i}(s)}W^n\left(s,\cdot,\mu^{n,i}(s)\right)\right\|^2\Big]ds\\\notag
&\quad+\frac{T}{4c_Bn^2}.
\end{align}
Combining the last bound with \eqref{COHEN555} we get
\begin{align}\notag
&\E\left[\left|W^n\left(t,X^{n,i}(t),\mu^{n,i}(t)\right)\right|^2\right]\\\notag
&\le \E\left[\left|W^n\left(t,X^{n,i}(t),\mu^{n,i}(t)\right)\right|^2\right]+\frac{1}{2c_B}\int_{t}^T\E\left[\left\|\Delta_{X^{n,i}(s)}W^n\left(s,\cdot,\mu^{n,i}(s)\right)\right\|^2\right]ds\\\notag
&\le 3C_1^2c_B\E\Big[\int_t^T
\left|W^n\left(s,X^{n,i}(s),\mu^{n,i}(s)\right)\right|^2ds\Big]
+\frac{T}{4c_Bn^2}.
\end{align}
Gr{\"o}nwall's inequality implies that 
\begin{align}\notag
&\sup_{t\in[0,T]}\E\left[\left|W^n\left(t,X^{n,i}(t),\mu^{n,i}(t)\right)\right|^2\right]\le Cn^{-2},
\end{align}
and therefore
\begin{align}\notag
&\frac{1}{2c_B}\int_{0}^T\E\left[\left\|\Delta_{X^{n,i}(s)}W^n\left(s,\cdot,\mu^{n,i}(s)\right)\right\|^2\right]ds\\\notag
&\qquad\le 3C_1^2c_B\E\Big[\int_0^T
\left|W^n\left(s,X^{n,i}(s),\mu^{n,i}(s)\right)\right|^2ds\Big]
+\frac{T}{4c_Bn^2}\le Cn^{-2}.
\end{align}
The last two bounds establish \eqref{225}.

\qed

\subsection{Fluctuations of the empirical measure}\label{sec25}
In this section we provide the dynamics of the process $\sqrt{n}(\mu^n-\mu)$ for some appropriate $\mu:[0,T]\to\calP([d])$, where
\begin{align}\label{236aaa}
\mu^{n}(t):=\frac{1}{n}\sum_{j\in[n]}e_{X^{n,j}(t)},\qquad t\in[0,T],
\end{align}
is the {\it Nash equilibrium empirical measure} of the states of all the players. That is, its generator is
\begin{align}\notag
\calL^n_tg(\eta)=\sum_{x,y\in[d]}n\eta_xa^*_y\left(x,\tfrac{n}{n-1}\eta-\tfrac{1}{n-1}e_x,\Delta_xV^{n,i}(t,\cdot,\tfrac{n}{n-1}\eta-\tfrac{1}{n-1}e_x)\right)\left(g(\eta+\tfrac{1}{n}e_{xy})-g(\eta)\right),
\end{align}
where $g:\calP^n([d])\to\R$. 
Throughout this subsection we assume that the following limit exists 
$\psi_0:=\lim_{n\to\iy}\sqrt{n}(\mu^n(0)-\mu(0))$. For example, for $\mu(0)=(p_x:x\in[d])\in\calP([d])$ and for every $n\in\N$ and $i\in[n]$, the random variables $X^{n,i}(0)$ are i.i.d.~with distribution $\mu(0)$, $\psi_0$ has a multivariate normal distribution with mean $0$ and a $d\times d$ covariance matrix, whose $ij$-th component is $p_i(\delta_{ij}-p_j)$, where $\delta_{ij}$ is the Kronecker delta. 
To establish the fluctuations, we start by coupling two jump processes whose transition rates are driven by the functions $V^{n}$ and $U$. They serve us to show that the difference between the empirical measures driven by these transition rates is of order ${\cal O}(n^{-1})$ and therefore we can restrict the fluctuations analysis to the dynamics that are driven by $U$.

\begin{construction}\label{const_21} Let $\{X^{n,j}(0):j\in[n]\}$ be a collection of independent identically distributed (i.i.d.) random variables taking values in $[d]$ and set $Y^{n,j}(0)=X^{n,j}(0)$, $j\in[n]$. Fix $i\in[n]$ and consider a Cox process $Q^{n,i}$ with rate 
\begin{align}\notag
a(t):=\max&\left\{\sum_{x\in[d]}a^X_x(t),\sum_{x\in[d]}a^Y_x(t)\right\},
\qquad t\in[0,T],
\end{align}
where 
\begin{align}\notag
a^X_x(t)&=\max\{0, a^*_x(X^{n,i}(t),\mu^{n,i}(t),\Delta_{X^{n,i}(t)}V^{n}(t,\cdot,\mu^{n,i}(t))\},\\\notag
a^Y_x(t)&=\max\{0, a^*_x( X^{n,i}(t),\nu^{n,i}(t),\Delta_{X^{n,i}(t)}U(t,\cdot,\nu^{n,i}(t))\},
\end{align}for $x\in[d]$, $t\in[0,T]$, and 
\begin{align}\notag
\mu^{n,i}(t):=\frac{1}{n-1}\sum_{j,j\ne i}e_{X^{n,j}(t)},\quad\nu^{n,i}(t):=\frac{1}{n-1}\sum_{j,j\ne i}e_{Y^{n,j}(t)},\quad (t,i)\in[0,T]\times[n].
\end{align}
We assume that the processes $\{Q^{n,i}\}_{i\in[n]}$ satisfy
\begin{align}\notag
\Big\langle Q^{n,i}(\cdot)-\int_0^\cdot a^{n,i}(s)ds, Q^{n,j}(\cdot)-\int_0^\cdot a^{n,j}(s)ds\Big\rangle(t)=\delta_{ij}\int_0^t a^{n,i}(s)ds,\qquad t\in[0,T],
\end{align}
where $\delta_{ij}=1$ if $i=j$ and $0$ otherwise. 

Let
\begin{align}\notag
\tau^{n,i}:=\inf\{t\ge 0: X^{n,i}(t)\ne Y^{n,i}(t)\}\wedge T,
\end{align}
with the convention that $\inf\emptyset=\iy$. 
On the time interval $[0,\tau^{n,i}]$, whenever $Q^{n,i}$ jumps, the 2-dimensional process
$\left(X^{n,i},Y^{n,i}\right)$ jumps to state 
\begin{align}\notag
 \begin{cases*}
     (x,x), & \text{with probability} $\quad\min\{a^X_x(t),a^Y_x(t)\}/a(t)$, \\
      \left(x,X^{n,i}(t)\right), & \text{with probability} $\quad\left(a^X_x(t)-\min\{a^X_x(t),a^Y_x(t)\}\right)/a(t)$, \\
     \left(X^{n,i}(t),x\right), & \text{with probability} $\quad\left(a^Y_x(t)-\min\{a^X_x(t),a^Y_x(t)\}\right)/a(t)$,
     \end{cases*}
\end{align}
Notice that in the last two cases at most one of the processes $X^{n,i}$ and $Y^{n,i}$ jumps.

On the time interval $[\tau^{n,i},T]$ the processes $X^{n,i}$ and $Y^{n,i}$ move independently according to the transition rates $(a^X_x(t):x\in[d])$ and $(a^Y_x(t):x\in[d])$, respectively.
\end{construction}
One can verify by induction over the jumps that the processes $Q^{n,i},X^{n,i}$, and $Y^{n,i}$, $i\in[n]$, are well-defined. 
The next proposition provides an approximation of order ${\cal O}(n^{-1})$ to the empirical measure in the $n$-player game. We also consider the empirical measures of the states of all the players $\mu^n$, which is given in \eqref{236aaa} and 
\begin{align}\label{241bb}
\nu^{n}(t):=\frac{1}{n}\sum_{j\in[n]}e_{Y^{n,j}(t)}.
\end{align}
That is, its generator is 
\begin{align}\notag
\calL^n_tg(\eta)=\sum_{x,y\in[d]}n\eta_xa^*_y\left(x,\tfrac{n}{n-1}\eta-\tfrac{1}{n-1}e_x,\Delta_xU(t,\cdot,\tfrac{n}{n-1}\eta-\tfrac{1}{n-1}e_x)\right)\left(g(\eta+\tfrac{1}{n}e_{xy})-g(\eta)\right)
\end{align}
where $g:\calP^n([d])\to\R$. 

\begin{proposition}\label{prop_21}
Under Assumptions \ref{assumption21} and \ref{assumption22}, there exists $C>0$ such that for every 
$n\ge 2$,
\begin{align}\label{240}
\E\Big[\sup_{s\in[0,T]}\left\|\mu^{n,i}(s)-\nu^{n,i}(s)\right\|\Big]\le 
2\E\Big[\sup_{s\in[0,T]}\left|X^{n,i}(s)-Y^{n,i}(s)\right|\Big]\le  \frac{C}{n},
\end{align}
and as a consequence
\begin{align}\label{241}
\E\Big[\sup_{s\in[0,T]}\left\|\mu^{n}(s)-\nu^{n}(s)\right\|\Big]\le \frac{C+2}{n}
\end{align}
and 
\begin{align}\label{241b}
\sqrt{n}\sup_{[0,T]}\|\mu^n-\nu^n\|\text{  converges in probability to $0$}.
\end{align}
\end{proposition}
{\bf Proof.}  The proof of \eqref{240} is the most demanding part, so we start with the other two. Inequality \eqref{241} merely follows from it, since 
$$\sup_{s\in[0,T]}\Big[\left\|\mu^{n}(s)-\mu^{n,i}(s)\right\|+\left\|\nu^{n}(s)-\nu^{n,i}(s)\right\|\Big]\le\frac{2}{n},$$ 
and \eqref{241b} follows from \eqref{241} by Markov's inequality.

We now turn to proving \eqref{240}. 
%
The first inequality follows since for every $(x_j:j\in[n]), (y_j:j\in[n])\in[d]^n$,
\begin{align}\notag
\Big\|\frac{1}{n}\sum_{j\in[n]}e_{x_j}-\frac{1}{n}\sum_{j\in[n]}e_{y_j}\Big\|\le \frac{2}{n}\sum_{j\in[n]}|x_j-y_j|
\end{align}
and since the processes $\{X^{n,j}-Y^{n,j}:j\in[n]\}$ are identically distributed. The inequality above follows since the Euclidean norm is bounded by the $l_1$ norm for which the inequality is straightforward. Therefore, we now turn to proving the second inequality. 
From Construction \ref{const_21} and the inequality $1-e^{-x}\le x$, it follows that for every $t\in[0,T]$ 
\begin{align}\notag
\E\Big[\sup_{s\in[0,t]}\left|X^{n,i}(s)-Y^{n,i}(s)\right|\Big]
&\le(d-1)\PP\Big(\tau^{n,i}\le t\Big)\\\notag
&= (d-1)\left(1-\E\Big[e^{-\sum_{z\in[d]}\int_{0}^t\left|a^X_z(s)-a^Y_z(s)\right|ds}\Big]\right)
\\\notag
&\le(d-1)\sum_{z\in[d]}\int_{0}^t\E\left[\left|a^X_z(s)-a^Y_z(s)\right|\right]ds.
\end{align}
Now, the Lipschitz-continuity of $a^*$ given in \eqref{225-a} implies that
\begin{align}\label{244}
&\E\left[\sup_{s\in[0,t]}\left|X^{n,i}(s)-Y^{n,i}(s)\right|\right]\\\notag
&\quad\le c_Ld(d-1)\left(\int_{0}^t\E\left[\left\|\mu^{n,i}(s)-\nu^{n,i}(s)\right\|\right]ds\right.\\\notag
&\qquad\qquad\qquad\qquad+
\int_{0}^t\E\left[\left\|
\Delta_{X^{n,i}(s)}V^{n}(s,\cdot,\mu^{n,i}(s))-\Delta_{X^{n,i}(s)}U(s,\cdot,\mu^{n,i}(s))\right\|\right]ds
\\\notag
&\qquad\qquad\qquad\qquad\left.+\int_{0}^t\E\left[\left\|
\Delta_{X^{n,i}(s)}U(s,\cdot,\mu^{n,i}(s))-\Delta_{X^{n,i}(s)}U(s,\cdot,\nu^{n,i}(s))\right\|\right]ds\right).
\end{align}
By the definition of the operator $\Delta_y$ and since $U$ is uniformly Lipschitz, we get that the last integral is bounded above by 
\begin{align}\notag
2\sqrt{d}c_L
\int_{0}^t\E\left[\left\|\mu^{n,i}(s)-\nu^{n,i}(s)\right\|\right]ds
&\le
\frac{2\sqrt{d}c_L}{n-1}\int_{0}^t\sum_{j,j\ne i}\E\left[\left|X^{n,j}(s)-Y^{n,j}(s)\right|\right]ds
\\\notag
&=2\sqrt{d}c_L\int_{0}^t\E\left[\left|X^{n,i}(s)-Y^{n,i}(s)\right|\right]ds,
\end{align}
where the equality follows by the symmetry of the players. The last bound applied to \eqref{244} together with \eqref{225} and Gr{\"o}nwall's inequality imply that there is a constant $C>0$ such that for every $t\in[0,T]$ and $n\ge 2$, 
\begin{align}\notag
\E\Big[\sup_{s\in[0,T]}\left|X^{n,i}(s)-Y^{n,i}(s)\right|\Big]\le \frac{C}{n}.
\end{align}

\qed

The next theorem provides the fluctuations of the process $\mu^n$. 
In the light of \eqref{241b}, our proof will focus on the fluctuations of $\nu^n$. Therefore, we set 
the $\calQ^{d\times d}$-valued function 
$\al^*$ by
\begin{align}\notag
\al^{*}_{xy}(s,\eta)&=a^*_y(x,\eta,\Delta_x U(s,\cdot,\eta)),\qquad x,y\in[d],
\end{align}
where pay attention that the operator $\Delta_x$ acts on $U$ and not on $V^{n}$. 
Also, the following notation is necessary for the statement of the theorem. For every $\calQ^{d\times d}$-valued function $\eta\mapsto\al(\eta)$, its gradient is a $d\times d$ matrix, whose $xy$ component is the vector $\nabla_\eta\al_{xy}(\eta)$. For any vector $c\in\R^d$, the product $c\otimes \nabla_\eta\al(\eta)$ is a $d\times d$ matrix, whose $xy$ component is the inner product $c\cdot\nabla_\eta\al_{xy}(\eta)$.  
Finally, in order to establish the fluctuations we require convergence of the initial state of the system and 
some regularity of the optimal control $a^*$.
\begin{assumption}\label{assumption23}
Modifying $c_L$ from Assumption \ref{assumption22} if necessary, we assume that for every $x,y\in[d]$, $\eta,\eta'\in\calP([d])$, and $p,p'\in[-2c_U,2c_U]^d$, one has
\begin{align}\notag
\|\nabla_\eta a^*_y(x,\eta,p)-\nabla_\eta a^*_y(x,\eta',p')\|+\|\nabla_p a^*_y(x,\eta,p)-\nabla_p a^*_y(x,\eta',p')\|&\le c_M\left(\|\eta-\eta'\|+\|p-p'\|\right).
 \end{align} 
\end{assumption}
The bound for the first part on the left-hand side of the above follows easily if for example the running cost can be expressed as $f(x,\eta,a)=f_1(x,a)+f_2(x,\eta)$ for some functions $f_1$ and $f_2$, in which case $a^*$ is independent of $\eta$. The optimal control in this case, $a^*(x,\eta,\Delta_x U(t,\cdot,\eta))$ depends on the empirical measure only through the function $U$. A sufficient condition for bounding the second term by the right-hand side (r.h.s.) is given by  in Assumption \ref{assumption1300}.5, for more details see its proceding paragraph.

\begin{theorem}\label{thm_22}
Let 
$(\mu(t))_{t\in[0,T]}$ be given by $\tfrac{d}{dt}\mu(t)^\top=\mu(t)^\top\al^*(t,\mu(t))$, with the given initial condition $\mu(0)=\mu_0\in\calP([d])$. Under Assumptions \ref{assumption21} and \ref{assumption22}, 
the process $\sqrt{n}(\mu^n-\mu)$ is stochastically bounded. That is, 
\begin{align}\label{250-z}
\lim_{k\to\iy}\sup_n\PP\left(\sup_{t\in[0,T]}\sqrt{n}\|\mu^n(t)-\mu(t)\|>k\right)=0.
\end{align}
Moreover, if in addition Assumption \ref{assumption23} holds, then  
$\sqrt{n}(\mu^n(\cdot)-\mu(\cdot))\To\psi(\cdot)$, where $\psi$ uniquely solves 
\begin{align}\label{250}
d\psi(t)&=\left[(\al^*(t,\mu(t)))^\top\psi(t)+\left(\psi(t)\otimes\nabla_\eta\al^*(t,\mu(t))\right)^\top\mu(t)\right]dt+\Sigma(t)dB(t),
\end{align}
on the interval $[0,T]$, with the initial condition $\psi(0)=\psi_0$.  
The process $B$ is a standard $d$-dimensional Brownian motion and $\Sigma:[0,T]\to\R^{d\times d}$ is given by
\begin{align}\notag
(\Sigma^2)_{xy}(t)&=
-\mu_y(t)\al^*_{yx}(t,\mu(t))-\mu_x(t)\al^*_{xy}(t,\mu(t)) ,\qquad x\ne y\\\notag
(\Sigma^2)_{xx}(t)&=\sum_{z,z\ne x}\mu_z(t)\al^{*}_{zx}(t,\mu(t))+\mu_x(t)\sum_{z,z\ne x}\al^{*}_{xz}(t,\mu(t)).
\end{align}
%
\end{theorem}
{\bf Proof.}  We start with the case that Assumption 
 \ref{assumption23} holds. The more general case is treated later. The stochastic differential equation admits a unique solution since the $dt$ component is linear in $\psi$ and the terms of $\al^*$ and $\nabla_\eta\al^*$ are bounded. We now turn to showing the convergence. 
First, note that by \eqref{241b} it is sufficient to show that $$\psi^n:=\sqrt{n}(\nu^n-\mu)\To\psi.$$
For every $x\in[d]$, let $A^n_x$ and $S^n_x$ be the arrival and departure processes, respectively, associated with state $x$. That is, $A^n_x(t)$ (resp., $S^n_x(t)$) counts how many times players moved into (from) state $x$ during the time interval $[0,t]$, so that $\sum_{x\in[d]}A^n_x(t)=\sum_{x\in[d]}S^n_x(t)$ is the total number of jumps of the process $\left(X^{n,i}:i\in[n]\right)$ during $[0,t]$. The rates of transition of $A^n_x(t)$ and $S^n_x(t)$ are, respectively, $n\la^n_x(t)$ and $n\sigma^n_x(t)$, where,
\begin{align}\notag
\la^n_x(t)&:=\sum_{y,y\ne x}\nu^n_y(t)\al^{*}_{yx}(t,\nu^{n,\sharp,y}(t)),\quad
\sigma^n_x(t):=\nu^n_x(t)\sum_{y,y\ne x}\al^{*}_{xy}(t,\nu^{n,\sharp,x}(t)),\quad t\in[0,T],
\end{align}
where for every $(x,\eta)\in[d]\times \calP^n([d])$ we set the $\calP([d])$ element,
\begin{align}\label{248aa}
\eta^{\sharp,x}:=
\begin{cases}
\frac{n}{n-1}\eta-\frac{1}{n-1}e_x, &\text{when $\eta_x>0$},\\
\eta, &\text{when $\eta_x=0$}.
\end{cases}
\end{align}
The term $\nu^{n,\sharp,x}(t)$ stands for the empirical distribution of the rest of the players from the point of view of a player whose state is $x$. 
Also, for every $x\in[d]$, set the martingales
\begin{align}\notag
M^n_x(t)=\frac{1}{\sqrt{n}}\Big(A^n_x(t)-\int_0^tn\lambda^n_x(s)ds \Big),\qquad
N^n_x(t)=\frac{1}{\sqrt{n}}\Big(S^n_x(t)-\int_0^tn\sigma^n_x(s)ds \Big),
\end{align}
where $t\in[0,T]$ and denote $M^n(t)=\left(M^n_x(t):x\in[d]\right)$, $t\in[0,T]$.

We start with the fluid scale and show that 
\begin{align}\label{290}
\sup_{t\in[0,T]}\|\nu^n(t)-\mu(t)\| \text{ converges in probability to $0$}.
\end{align}
From the limit $\sqrt{n}(\mu^n(0)-\mu(0))\To0$ and \eqref{241b}, 
\begin{align}\label{290a}
\nu^n(0)\To\mu(0).
\end{align}
Simple algebraic manipulations yield that
\begin{align}\label{291}
\nu^n(t)-\mu(t)&=\nu^n(0)-\mu(0)+\frac{1}{\sqrt{n}}(M^n(t)-N^n(t))\\\notag
&\quad+\int_0^t\left[\nu^n(s)^\top\al^{*}(s,\nu^{n,\sharp}(s))-\mu(s)^\top\al^{*}(s,\mu(s))\right]ds,
\end{align}
where $\al^*_{xy}(s,\nu^{n,\sharp}(s))$ reads as $\al^*_{xy}(s,\nu^{n,\sharp,x}(s))$. 
By Assumption \ref{assumption22} and the bound 
\begin{align}\label{292}
\sup_{(t,x)\in[0,T]\times[d]}\|\nu^n(t)-\nu^{n,\sharp,x}(t)\|\le 2(n-1)^{-1},
\end{align}
we obtain that
\begin{align}\notag
\sup_{s\in[0,t]}\|\nu^n(s)-\mu(s)\|&\le \|\nu^n(0)-\mu(0)\|+\frac{1}{\sqrt{n}}\sup_{s\in[0,t]}\|M^n(s)-N^n(s)\|\\\notag&\quad+C\Big(\int_0^t\sup_{u\in[0,s]}\|\nu^n(u)-\mu(u)\|ds+\frac{1}{n}\Big), 
\end{align}
where in the above expression, and in the rest of the proof, $C$ refers to a finite positive constant that is independent of $t$ and $n$ and which can change from one line to the next. 
Now, Gr{\"o}nwall's inequality implies that 
\begin{align}\label{254-a}
\sup_{s\in[0,T]}\|\nu^n(s)-\mu(s)\|\le C\Big(\|\nu^n(0)-\mu(0)\|+\frac{1}{\sqrt{n}}\sup_{s\in[0,T]}\|M^n(s)-N^n(s)\|+\frac{1}{n}\Big). 
\end{align}
Since any pair of $\{A^n_x,S^n_y:x,y\in[d]\}$ is orthogonal, see \eqref{202}, we have for every $x\in[d]$, 
\begin{align}\label{254}
\langle M^n_x-N^n_x\rangle(t)&=\frac{1}{n}\Big\langle A^n_x(\cdot)-\int_0^\cdot n\lambda^n_x(s)ds \Big\rangle(t)+\frac{1}{n}\Big\langle S^n_x(\cdot)-\int_0^\cdot n\sigma^n_x(s)ds \Big\rangle(t)\\\notag
&=
\int_0^t\la^n_x(s)ds+\int_0^t\sigma^n_x(s)ds.
\end{align}
Since $\la^n_x$ and $\sigma^n_x$ are uniformly bounded we get that 
$$\sup_{t\in[0,T]}\left\langle\frac{1}{\sqrt{n}}(M^n_x-N^n_X)\right\rangle(t)\le \frac{C}{n}.$$
Therefore, from Burkholder--Davis--Gundy inequality (\cite[Theorem 48]{Protter2004}), \eqref{254-a}, \eqref{290a}, and Markov inequality we get \eqref{290}.

We now study the diffusion scale. From \eqref{291}
\begin{align}\label{252}
\psi^n(t)-\psi^n(0)
&=M^n(t)-N^n(t)+L^n(t)+K^n(t),
\end{align}
where
\begin{align}\label{253a}
L^n(t)&:=\int_0^t\psi^n(s)^\top\al^{*}(s,\nu^{n,\sharp}(s))ds,\\\label{253b}
K^n(t)&:=\int_0^t\sqrt{n}\mu(s)^\top\left[\al^{*}(s,\nu^{n,\sharp}(s))-\al^{*}(s,\mu(s))\right]ds.
\end{align}

To attain convergence, we use tightness arguments.\footnote{The following definition is borrowed from \cite[Definition VI.3.25]{Jacod1987}. 
A sequence of stochastic processes with sample paths in $\calD([0,T],\R^k)$, $k\in\N$, is said to be $\calC$-tight if it is tight and every subsequential limit has continuous sample paths with probability 1. 
} 
In the next few paragraphs we will also use the following characterization of $\calC$-tightness for processes with sample paths in $\calD([0,T],\R^k)$,  see \cite[Proposition VI.3.26]{Jacod1987}: The sequence  $\{F^n\}_n$ is $\calC$-tight if and only if the sequence $\{\sup_{0\le t\le T}\|F^n(t)\|\}_n$ is tight and for every $\eps > 0$ and $\gamma > 0$ there exist $N_0$ and $\theta > 0$ such that
for every $n \ge N_0$, 
$$P\Big(\sup_{0\le s\le t\le (s+\theta)\wedge T}\|F^n(t)-F^n(s)\|>\gamma\Big)<\eps.$$
Next, we argue the $\calC$-tightness 
 of $\{(\psi^n,M^n-N^n,L^n+K^n)\}$ in $\calD([0,T],\R^{3d})$. 
From \cite[Ch.~VI, Corollary 3.33]{Jacod1987} and \eqref{252} it is sufficient to show separately the $\calC$-tightness of $\{M^n-N^n\}$ and $\{L^n+K^n\}$. Observe that merely tightness of the last two sequences does not imply the tightness of their sum or their joint distribution. Hence, we appeal to $\calC$-tightness. 
We start with $\{M^n-N^n\}$. 
Notice that for any $x\ne y$,
\begin{align}\notag
\langle M^n_x-N^n_x,M^n_y-N^n_y\rangle(t)&=-\langle M^n_x,N^n_y\rangle(t)-\langle M^n_y,N^n_x\rangle(t)\\\notag
&=-\int_0^t[\nu^n_y(s)\al^*_{yx}(s,\nu^{n,\sharp,y}(s))+\nu^n_x(s)\al^*_{xy}(s,\nu^{n,\sharp,x}(s))]ds.
\end{align}
The Martingale central limit theorem (\cite[Ch.~7, Theorem 1.4]{Ethier1986}), which holds due to \eqref{290}, \eqref{254}, \eqref{292}, and the last equation implies that 
\begin{align}\notag
(M-N)(\cdot):=\lim_{n\to\iy}(M^n-N^n)(\cdot)=\int_0^\cdot\Sigma(t)dB(t). 
\end{align}
Therefore, $\{M^n-N^n\}$ is $\calC$-tight.

We now turn to the $\calC$-tightness of $\{L^n+K^n\}$. First, Assumption \ref{assumption23}, \eqref{225-c}, 
Lemma \ref{lem_a1}, and the boundedness of $\nabla_\eta\al^*$, as a continuous function, imply, 
\begin{align}\label{255}
K^n(t)&=\int_0^t\mu(s)^\top\left[\psi^n(s)^\top\otimes\nabla_\eta\al^{*}(s,\mu(s))\right]ds\\\notag
&\quad+{\cal O}\Big(\sup_{(s,x)\in[0,T]\times[d]}\|\psi^n(s)\|\|\nu^{n,\sharp,x}(s)-\mu(s)\|\Big).
\end{align}
Next, we wish to bound $\|\psi^n(t)\|$. From the definition of $L^n$, the last representation of $K^n$, and again the boundedness of $\al^*$ and $\nabla_\eta\al^*$ it follows that $\|L^n(t)+K^n(t)\|\le C\int_0^t\|\psi^n(s)\|ds.$
Going back to \eqref{252} and using Gr{\"o}nwall's inequality, we get the following estimation
\begin{align}\notag
\|\psi^n(t)\|\le C\|M^n(t)-N^n(t)\|,\qquad t\in[0,T].
\end{align}
The tightness of $\{M^n-N^n\}$ implies that $$\lim_{k\to\iy}\limsup_n\PP\left(\sup_{t\in[0,T]}\|M^n(t)-N^n(t)\|\ge k\right)=0.$$ Then, together with \eqref{253a} and \eqref{255}, the boundedness of $\al^*$ and the elements of $\nabla_\eta\al^*$, and \eqref{292}, we get that $\{L^n+K^n\}$ is $\calC$-tight in $\calD([0,T],\R^d)$.

The last step is to show that any weak limit $\psi$ of $\psi^n$ satisfies \eqref{250}.
Since $\{(\psi^n,M^n-N^n,L^n+K^n)\}$ is tight it has a convergence subsequence, which we relabel as $\{n\}$, with limit $\{(\psi,M-N,L+K)\}$. Then, \eqref{252} yields
$$
\psi(t)-\psi(0)=(M-N)(t)+(L+K)(t),\qquad t\in[0,T].$$
From \eqref{290}, we get that for every $x\in[d]$, $(\nu^n,\nu^{n,\sharp,x})$ converges in probability to $(\mu,\mu)$. Therefore,
\begin{align}\notag
(L+K)(t)&=\int_0^t\Big[\psi(s)^\top\al^{*}(s,\mu(s))+\mu(s)^\top\left(\psi(s)^\top\otimes\nabla_\eta\al^{*}(s,\mu(s))\right)\Big]ds.
\end{align}

In case that Assumption \ref{assumption23} does not hold, using merely the Lipschitz-continuity of $a^*$, the approximation \eqref{255} to \eqref{253b} is replaced by the bound
\begin{align}\notag
\|K^n(t)\|\le C\int_0^t\|\psi^n(s)\|ds,
 \end{align} 
 for some constant $C>0$ independent of $n$ and $t$. The boundedness of $a^*$ implies the same bound for $L^n(t)$. Now from \eqref{252}, Gr{\"o}nwall's inequality, and the $\calC$-tightness of $\{M^n-N^n\}$, we obtain \eqref{250-z}.

\qed

\section{Sufficient conditions for Theorems \ref{thm_21} and \ref{thm_22}
}\label{sec1300}
\beginsec

In this section we show that a variant of the conditions imposed in \cite{Gomes2013} are sufficient for Assumptions \ref{assumption21}--\ref{assumption23} to hold. Throughout this process we explicitly show how to derive the master equation from the MFG. 
\begin{assumption}\label{assumption1300}
\begin{enumerate}
\item (control set) The minimal rate of transition allowed from one state to another is positive.
\item (separation) The function $f$ can be expressed as $f(x,\eta,a)=f_1(x,a)+f_2(x,\eta)$, $(x,\eta,a)\in[d]\times\calP([d])\times\R^d$.
\item (Lipschitz-continuity)
 The functions $f_2$ and $g$ are differentiable w.r.t.~$\eta$ and there exists a constant $c_L$ such that for every $x\in[d]$, and $\eta,\eta'\in\calP([d])$,
\begin{align}\notag
&\|\nabla_\eta g(x,\eta)-\nabla_\eta g(x,\eta')\|+\|\nabla_\eta f_2(x,\eta)-\nabla_\eta f_2(x,\eta')\|\le c_L\|\eta-\eta'\|.
\end{align}
\item (monotonicity) For every $\eta,\eta'\in\calP([d])$, 
\begin{align}\notag
\sum_{x\in[d]}(\eta_x-\eta_x')(g(x,\eta)-g(x,\eta'))&\ge 0,\\\notag
\sum_{x\in[d]}(\eta_x-\eta_x')(f_2(x,\eta)-f_2(x,\eta'))&\ge 0.
\end{align}
\end{enumerate}
\noi From the representation of $f$ it follows that $H(x,\eta,p)=H_1(x,p)+f_2(x,\eta)$, where
$\\$$H_1(x,p):=\inf_{a} \{f_1(x,a)+\sum_{y,y\ne x}a_yp_y\}$, and the infimum is taken over the controls that satisfy property 1 above.
\begin{enumerate}[resume]
\item (concavity) The function $H_1$ is twice continuously differentiable. Moreover, for every $M>0$ there exists a parameter $c_M>0$ such that for every $x,y,z\in[d]$ and $p\in[-M,M]^d$,
 \begin{align}\notag
\partial^2_{p_x,p_y}H_1(z,p)\le - c_M,
\end{align}
where $\partial^2_{p_x,p_y}$ stands for the second partial derivative. Moreover, for every $x,y,z\in[d]$, $p\mapsto \partial^2_{p_x,p_y}H_1(z,p)$ is Lipschitz continuous on $[-M,M]^d$. 
\end{enumerate}
\end{assumption}
Notice that Assumption \ref{assumption1300}, which is in force in \cite{CDLL2015}, implies the conditions given in Lemmas \ref{lem21} and \ref{lem22}. 
The conditions here are stronger than the ones given in Lemma \ref{lem21} and Remark \ref{rem21} as we consider now a lower bound for the control set and a more regular structure. 
By the regularity of $H_1$,
$a^*_x(y,p)=\partial_{p_x}H_1(y,p)$ and $\partial_{p_z}a^*_x(y,p)=\partial^2_{p_z,p_x}H_1(y,p)$. By the regularity and the monotonicity of $f_2$ and $g$, for every $p\in\calM([d])=\{\eta\in\R^d:\sum_{x\in[d]}\eta_x=0\}$,
\begin{align}\label{ASAF010}
&\sum_{x\in[d]}p_x\nabla_\eta g(x,\mu(T))\cdot p\ge 0,\\\label{ASAF011}
&\sum_{x\in[d]}p_x\nabla_\eta f_2(x,\mu(t))\cdot p\ge 0.
\end{align}
Notice also that by the separation condition, for any $x,y\in[d]$, $\eta,\eta'\in\calP([d])$, and $p\in\R$, $a^*_x(y,\eta,p)=a^*_x(y,\eta',p)$. Therefore, we abuse notation and in this section we use the notation $a^*_x(y,p)$. As a result Assumption \ref{assumption23} holds trivially.

The only seemingly missing parts from Assumptions \ref{assumption21} and \ref{assumption22} are that the master equation \eqref{215} has a unique classical solution with a Lipschitz gradient. 
The next theorem states that these conditions are guaranteed. 
Before stating the theorem, we provide an example that satisfies the conditions of Assumption \ref{assumption1300}.
\begin{example}\label{exmaple1300}
Fix a  control set $[l,L]$. Consider twice continuously differentiable and convex functions $b_1,b_2:[0,1]\to\R$. Set $g(x,\eta)=b_1(\eta_x)$ and $f_2(\eta)=b_2(\eta_x)$. Also, set a function $b_3:[d]\to\R$ and constants $\{c_y\}_{y\in[d]}$, $\{q_y\}_{y\in[d]}$, where $q_y\in(l,L)$, $x\in[d]$.  
Consider the function $f_1(x,a)=b_1(x)+\sum_{y,y\ne x}c_y(a_y-q_y)^2$. Then for sufficiently large $\{c_x\}_x$, the fifth condition holds. 
\end{example}


\begin{theorem}\label{thm1300}
Under Assumption \ref{assumption1300} the master equation \eqref{215} admits a unique classical solution $U$ and its gradient with respect to the measure argument, $\nabla_\eta U$, is continuous over $[0,T]\times[d]\times\calP([d])$. Moreover, 
there exists $c_L>0$ such that for every $t\in[0,T]$, $x\in[d]$, and $\eta,\eta'\in\calP([d])$,
\begin{align}\label{225-ccc}
\|\nabla_\eta U(t,x,\eta)-\nabla_\eta U(t,x,\eta')\|\le c_L\|\eta-\eta'\|.
\end{align}
\end{theorem}

The proof of the theorem is done in several steps and is given after a few preliminary lemmas and propositions. In the rest of this section we assume that Assumption \ref{assumption1300} is in force. The proof is inspired by \cite{CDLL2015} 
and its idea is as follows. We first consider the forward-backward system given in \eqref{224united}, and use its solution to define a function $U$. We show that $U$ is differentiable w.r.t.~the measure component and satisfies the master equation \eqref{215}. Finally, we show its regularity.
The function $u(t,x)$ stands for $V(t,x,\mu(t))$ from \eqref{224united}. 
Also, set $(t_0,\mu_0)\in[0,T]\times\calP([d])$ and consider the forward-backward system:
\begin{align}
\begin{cases}\label{1300}
\frac{d}{dt} \mu_x(t)=\sum_{y\in[d]}\mu_y(t)a^*_x(y,\Delta_y u(t,\cdot)),\\
-\frac{d}{dt} u(t,x)=H_1(x,\Delta_xu(t,\cdot))+f_2(x,\mu(t)),\\
u(T,x)=g(x,\mu(T)),\quad \mu(t_0)=\mu_0,
\end{cases}
\end{align}
where $\mu:[t_0,T]\to\calP([d])$ and $u:[t_0,T]\times[d]\to\R$. The $x$-th component of $\mu(t)$ is denoted by $\mu_x(t)$, $x\in[d]$. This system has a unique classical solution; see Lemma \ref{lem22}. 
Every function $\phi:[d]\to\R$ can be identified with a $d$-dimensional vector $\phi\in\R^d$. Hence $\|\phi(\cdot)\|$ should be understood as the Euclidean norm of its vector presentation. 
The next lemma provides a sensitivity result for the last system.
\begin{lemma}\label{lem1301}
Fix  $0\le t_0\le T$ and $\mu_0,\hat\mu_0\in\calP([d])$. Let $(\mu,u)$ and $(\hat\mu,\hat u)$ be two solutions to \eqref{1300} with the initial conditions $\mu(t_0)=\mu_0$ and $\hat\mu( t_0)=\hat\mu_0$. 
Then, there is a constant $C>0$ independent of $t_0$, $\mu_0$, and $\hat\mu_0$, such that 
\begin{align}\label{1300b}
&\sup_{t\in[t_0,T]}\left(\|\mu(t)-\hat\mu(t)\|
+\|u(t,\cdot)-\hat u(t,\cdot)\|\right)
\le C
\|\mu_0-\hat\mu_0\|.
\end{align}
\end{lemma}
{\bf Proof.} 
From \eqref{1300}, 
\begin{align}\notag
\frac{d}{dt}(\mu_x(t)-\hat\mu_x(t))&=\sum_{y\in[d]}(\mu_y(t)-\hat\mu_y(t))a^*_x(y,\Delta_yu(t,\cdot))
\\\notag
&\quad+\sum_{y\in[d]}\hat\mu_y(t)[a^*_x(y,\Delta_yu(t,\cdot))-a^*_x(y,\Delta_y\hat u(t,\cdot))].
\end{align}
From Assumption \ref{assumption22}.3 and Gr{\"o}nwall's and Jensen's inequalities, for every $t\in[t_0,T]$,
\begin{align}\label{1302}
\|\mu(t)-\hat\mu(t)\| &\le C\Big(\|\mu(t_0)-\hat\mu(t_0)\|+\int_{t_0}^t\sum_{y\in[d]}\hat\mu_y(s)\|\Delta_y u(s,\cdot)-\Delta_y \hat u(s,\cdot)\|ds\Big)
\\\notag
&\le C\|\mu(t_0)-\hat\mu(t_0)\|+C\Big(
\int_{t_0}^t\sum_{y\in[d]}\hat\mu_y(s)\|\Delta_y u(s,\cdot)-\Delta_y \hat u(s,\cdot)\|^2ds
\Big)^{1/2},
\end{align}
where here and in the rest of the proof, $C$ refers to a finite positive constant that is independent of $t_0$, $\hat t_0$, $\mu_0$, and $\hat\mu_0$ and which can change from one line to the next.

Now, by \cite[Equation (20)]{Gomes2013}, Assumptions \ref{assumption1300}.4 and  \ref{assumption1300}.5 imply that for every $t\in[t_0,T]$,
\begin{align}\notag
\sum_{y\in[d]}\hat\mu_y(t)\|\Delta_y u(t,\cdot)-\Delta_y \hat u(t,\cdot)\|^2
\le -\frac{d}{dt}\left[(\mu(t)-\hat\mu(t))\cdot(u(t,\cdot)-\hat u(t,\cdot))\right].
\end{align}
Integrating both sides and recalling Assumption \ref{assumption1300}.4, we obtain
\begin{align}\notag
\int_{t_0}^T\sum_{y\in[d]}\hat\mu_y(s)\|\Delta_y u(s,\cdot)-\Delta_y \hat u(s,\cdot)\|^2ds
\le \|\mu(t_0)-\hat\mu(t_0)\|\|u(t_0,\cdot)-\hat u(t_0,\cdot) \|,
\end{align}
Plugging this into \eqref{1302}, we get that
\begin{align}\notag
\sup_{t\in[t_0,T]}\|\mu(t)-\hat\mu(t)\|\le C\Big[\|\mu(t_0)-\hat\mu(t_0)\|+ (\|\mu(t_0)-\hat\mu(t_0)\|\|u(t_0,\cdot)-\hat u(t_0,\cdot)\| )^{1/2}\Big].
 \end{align} 
%
Using the Lipschitz continuity of $p\mapsto H_1(\cdot,p)$, $\eta\mapsto f_2(\cdot,\eta)$, and $\eta\mapsto g(\cdot,\eta)$ and Gr{\"o}nwall's inequality, we can derive from \eqref{1300} that 
\begin{align}\label{1306}
\sup_{t\in[t_0,T]}\|u(t,\cdot)-\hat u(t,\cdot))\|\le C\sup_{t\in[t_0,T]}\|\mu(t)-\hat \mu(t))\|,
 \end{align} 
which together with the previous bound implies \eqref{1300b}.

\qed

Since $H_1,f_2$, and $g$ are Lipschitz it follows from \eqref{1300} that $u$ is bounded. Using the continuity of $p\mapsto (a^*_x(y,p),\nabla_pa^*_x(y,p))$, which follows since $H_1$ is twice continuously differentiable, we get that $a^*_x(y,\Delta_yu(t,\cdot))$ and $\nabla_pa^*_x(y,\Delta_yu(t,\cdot))$ are bounded. We will make use of these properties more than once in the sequel.

We now introduce a {\it linearized forward-backward system}, $(\rho,w)$, which will be in use several times in the analysis. Recall the definition of $\calM([d])$ given before \eqref{ASAF010} and set measurable and bounded functions $k:[t_0,T]\to\calM([d])$, $l:[t_0,T]\to\R^d$, and $\zeta:[d]\to\R$. 
Also, set $\rho_0\in\calM([d])$, and consider a solution to \eqref{1300}, denoted by $(\mu,u)$, associated with the initial condition $\mu(t_0)=\mu_0\in\calP([d])$. Let $(\rho,w)$ be the solution of 
\begin{align}
\begin{cases}\label{1309}
\frac{d}{dt} \rho_x(t)=\sum_{y\in[d]}\rho_y(t)a^*_x(y,\Delta_y u(t,\cdot))+\sum_{y\in[d]}\mu_y(t)\nabla_pa^*_x(y,\Delta_yu(t,\cdot))\cdot\Delta_yw(t,\cdot)+k_x(t),\\
-\frac{d}{dt} w(t,x)=
\Delta_xw(t,\cdot)\cdot a^*(x,\Delta_xu(t,\cdot))
+\nabla_\eta f_2(x,\mu(t))\cdot \rho(t)+l_x(t),\\
w(T,x)=\nabla_\eta g(x,\mu(T))\cdot \rho(T)+\zeta(x),\quad \rho(t_0)=\rho_0,
\end{cases}
\end{align}
where $\rho:[t_0,T]\to\R^d$ and $w:[t_0,T]\times[d]\to\R$. The $x$-th components of $\rho(t)$, $k(t)$, and $l(t)$ are denoted by $\rho_x(t)$, $k_x(t)$, and $l_x(t)$, $x\in[d]$. 
\begin{proposition}\label{prop31}
The system \eqref{1309} has a unique classical solution. Moreover, there is a constant $C>0$, independent of $t_0, \rho_0, \zeta,k$, and $l$, such that
\begin{align}\label{ASAF004}
\sup_{t\in [t_0,T]}\left(\|\rho(t)\|+\|w(t,\cdot)\|\right)\le C\Big( \|\rho_0\|+\|\zeta\|+\sup_{t\in [t_0,T]}\left(\|k(t)\|+\|l(t)\|\right)\Big).
 \end{align} 
\end{proposition}
{\bf Proof.}  The proof appeals to 
the Leray--Schauder fixed point theorem: 
{\it let $\Phi$ be a continuous and compact mapping of a Banach Space $\mathbb{X}$ to itself such that the set $\mathcal{A}=\{\rho\in\mathbb{X}:\rho=\lambda \Phi(\rho)\;\;\text{for some}\;\;0\le \lambda\le 1\}$ is bounded. Then $\Phi$ has a fixed point.}

We now define a mapping $\Phi:\calC([t_0,T],\calM([d]))\to\calC([t_0,T],\calM([d]))$, where $\calC([0,T],\calM([d]))$ is the space of continuous functions, mapping $[t_0,T]\to\calM([d])$. 
Fix $\rho\in\calC([t_0,T],\calM([d]))$ and let $w$ be the classical solution of 
\begin{align}
\begin{cases}\label{ASAF020}
-\frac{d}{dt} w(t,x)=
\Delta_xw(t,\cdot)\cdot a^*(x,\Delta_xu(t,\cdot))
+\nabla_\eta f_2(x,\mu(t))\cdot \rho(t)+l_x(t),\\
w(T,x)=\nabla_\eta g(x,\mu(T))\cdot \rho(T)+\zeta(x).
\end{cases}
\end{align}
Given the solution $w$, set $\Phi(\rho)=\hat\rho$ to be the solution of 
\begin{align}
\begin{cases}\notag
\frac{d}{dt} \hat\rho_x(t)=\sum_{y\in[d]}\hat\rho_y(t)a^*_x(y,\Delta_y u(t,\cdot))+\sum_{y\in[d]}\mu_y(t)\nabla_pa^*_x(y,\Delta_yu(t,\cdot))\cdot\Delta_yw(t,\cdot)+k_x(t),\\
\hat \rho(t_0)=\rho_0.
\end{cases}
\end{align}
Notice that since $\sum_{x\in[d]}a^*_x(y,p)=0$ and by the choice of $k$ and $\rho_0$, we get that $\hat\rho\in\calC([t_0,T],\calM([d]))$. 
The mapping $\Phi$ is clearly continuous. In the rest of the proof we will demonstrate that $\mathcal{A}$ is bounded. It will then be clear also that \eqref{ASAF004} holds.
First let us show
\begin{align}\label{ASAF005}
\sup_{t\in [t_0,T]}\|w(t,\cdot)\|\le C\Big( \sup_{t\in[t_0,T]}\|\rho(t)\|+\|\zeta\|+\sup_{t\in [t_0,T]}\|l(t)\|\Big),
 \end{align} 
where, here and in the rest of the proof, $C$ refers to a generic positive constant, independent of $t_0, \rho_0, \zeta,k$, and $l$, and which can change from one line to the next. Fix $t\in[t_0,T]$. Integrating both sides of \eqref{ASAF020} and recalling that $a^*$, $\nabla_\eta f_2$, and $\nabla_\eta g$ are bounded, we get that for every $x\in[d]$,
\begin{align}\notag
w(t,x)&=\nabla_\eta g(x,\mu(T))\cdot\rho(T)+\zeta(x)\\\notag
&\quad+\int_t^T\Big[\Delta_xw(s,\cdot)\cdot a^*(x,\Delta_xu(s,\cdot))
+\nabla_\eta f_2(x,\mu(s))\cdot \rho(s)+l_x(s)\Big]ds\\\notag
&\le C\|\rho(T)\|+\|\zeta\|+C\int_t^T\Big[\|w(s,\cdot)\|
+\| \rho(s)\|+\|l(s)\|\Big]ds\end{align}
Therefore,
\begin{align}\notag
\|w(t,\cdot)\|\le C\Big(\|\zeta\|+\sup_{s\in[t_0,T]}\left(\| \rho(s)\|+\|l(s)\|\right)+\int_t^T\|w(s,\cdot)\|ds\Big).
\end{align}
Gr{\"o}nwall's inequality gives \eqref{ASAF005}.

Fix $\la\in[0,1]$. The identity $\rho=\la \Phi(\rho)$ implies that $(\rho,w)$ satisfies
\begin{align}
\begin{cases}\notag
\frac{d}{dt} \rho_x(t)=\sum_{y\in[d]}\rho_y(t)a^*_x(y,\Delta_y u(t,\cdot))\\
\quad\quad\quad\quad\quad+\la\Big(\sum_{y\in[d]}\mu_y(t)\nabla_pa^*_x(y,\Delta_yu(t,\cdot))\cdot\Delta_yw(t,\cdot)+k_x(t)\Big),\\
-\frac{d}{dt} w(t,x)=
\Delta_xw(t,\cdot)\cdot a^*(x,\Delta_xu(t,\cdot))
+\nabla_\eta f_2(x,\mu(t))\cdot \rho(t)+l_x(t),\\
w(T,x)=\nabla_\eta g(x,\mu(T))\cdot \rho(T)+\zeta(x),\quad \rho(t_0)=\la\rho_0,
\end{cases}
\end{align}
Therefore,
\begin{align}\notag
&\frac{d}{dt}\sum_{x\in[d]}w(t,x)\rho_x(t)\\\notag
&\quad=-\sum_{x\in[d]} \rho_x(t)\Delta_xw(t,\cdot)\cdot a^*(x,\Delta_xu(t,\cdot))-\sum_{x\in[d]}\rho_x(t)\nabla_\eta f_2(x,\mu(t))\cdot \rho(t)-\sum_{x\in[d]}\rho_x(t) l_x(t)\\\notag
&\quad+\la\sum_{x,y\in[d]}w(t,x)\mu_y(t)\nabla_pa^*_x(y,\Delta_yu(t,\cdot))\cdot\Delta_yw(t,\cdot)+\la\sum_{x\in[d]}w(t,x)k_x(t)\\\notag
&\qquad+\sum_{x,y\in[d]}w(t,x)\rho_y(t)a^*_x(y,\Delta_y u(t,\cdot))
.
\end{align}
The first and last terms on the r.h.s.~sum up to $0$. This is because $\sum_{y\in[d]}a^*_y(x,p)=0$.
Integrating both sides and rearranging, we get
\begin{align}\label{ASAF008}
&\la\int_{t_0}^T\sum_{x,y\in[d]}w(t,x)\mu_y(t)\nabla_pa^*_x(y,\Delta_yu(t,\cdot))\cdot\Delta_yw(t,\cdot)dt\\\notag
&\quad=\sum_{x\in[d]}\rho_x(T)\nabla_\eta g(x,\mu(T))\cdot\rho(T)+\zeta\cdot\rho(T)-\la w(t_0,\cdot)\cdot \rho_0
\\\notag
&\qquad+
\sum_{x\in[d]}\int_{t_0}^T\Big[\rho_x(t)\nabla_\eta f_2(x,\mu(t))\cdot \rho(t)+\rho_x(t)l_x(t)-\la w(t,x)k_x(t)\Big]dt\\\notag
&\quad\ge \zeta\cdot\rho(T)-\la w(t_0,\cdot)\cdot \rho_0+
\sum_{x\in[d]}\int_{t_0}^T\Big[\rho_x(t)l_x(t)-\la w(t,x)k_x(t)\Big]dt,
\end{align}
where the inequality follows from \eqref{ASAF010} and \eqref{ASAF011}.  
Using again $\sum_{y\in[d]}a^*_y(x,p)=0$, one may express the sum within the integral on the left-hand side as
\begin{align}\notag
-\sum_{y\in[d]}\mu_y(t)(\Delta_yw(t,\cdot))^\top A(y)\Delta_yw(t,\cdot),
\end{align}
where $A(y)$ is a $d\times d$ matrix, whose components are given by $A_{x,z}(y)=-\partial_{z}a^*_{x}(y,\Delta_y u(t,\cdot))=-\partial^2_{p_x,p_z}H(y,\Delta_y u(t,\cdot))$. Notice that $A(y)$ is symmetric and positive semidefinite. 
Thus, 
\begin{align}\label{ASAF009}
&\la\int_{t_0}^T\sum_{y\in[d]}\mu_y(t)(\Delta_yw(t,\cdot))^\top A(y)\Delta_yw(t,\cdot),
\\\notag
&\quad\le 
C\Big(\sup_{t\in[t_0,T]}\|\rho(t)\|\Big(\|\zeta\|+\sup_{t\in[t_0,T]}\|l(t)\|\Big)+\la\sup_{t\in[t_0,T]}\|w(t,\cdot)\|\Big(\|\rho_0\|+\sup_{t\in[t_0,T]}\|k(t)\|\Big)\Big).
\end{align}
Fix $\xi:[d]\to\R$ and $s\in(t_0,T]$ and consider the system
\begin{align}
\notag
-\frac{d}{dt} \hat w(t,x)=
\Delta_x\hat w(t,\cdot)\cdot a^*(x,\Delta_xu(t,\cdot)) \;\;\text{on}\;\;[t_0,s]  ,\qquad
w(s,x)=\xi(x).
\end{align}
The same way we obtained \eqref{ASAF005}, we may obtain the bound $\sup_{t\in[t_0,s]}\|\hat w(t,\cdot)\|\le C\|\xi\|$, where $C>0$ is independent of $\xi$ and $s$. 
Repeating the same steps as above, now for $\frac{d}{dt}\sum_{x\in[d]}\hat w(t,x)\rho_x(t)$, and using that $\mu(t)\in\calP([d])$, we get that
\begin{align}\notag
\rho(s)\cdot\xi(s)&=\la\hat w(t_0,\cdot)\cdot\rho_0+\la\sum_{x\in[d]}\int_{t_0}^s\hat w(t,x)k_x(t)dt\\\notag
&\quad-\la\int_{t_0}^s\sum_{x\in[d]}\mu_x(t)(\Delta_x\hat w(t,\cdot))^\top A(x)\Delta_x w(t,\cdot)dt\\\notag
&\le \la C\Big[\|\xi\|\|\rho_0\|+\|\xi\|\sup_{t\in[t_0,s]}\|k(t)\|\\\notag
&\qquad\qquad+\Big(\int_{t_0}^s\sum_{x\in[d]}\mu_x(t)(\Delta_x\hat w(t,\cdot))^\top A(x)\Delta_x \hat w(t,\cdot)dt\Big)^{1/2}\\\notag
&\qquad\qquad\qquad\times\Big(\int_{t_0}^s\sum_{x\in[d]}\mu_x(t)(\Delta_x w(t,\cdot))^\top A(x)\Delta_x w(t,\cdot)dt\Big)^{1/2}\Big]
.
\end{align}
From \eqref{ASAF009} and since the elements of $A(y)$ are bounded and $\sup_{t\in[t_0,s]}\|\hat w(t,\cdot)\|\le C\|\xi\|$, we get that the r.h.s.~of the above is bounded by 
\begin{align}\notag
&C\|\xi\|\Big(
\la\|\rho_0\|+\la\sup_{t\in[t_0,T]}\|k(t)\|+\la^{1/2}\sup_{t\in[t_0,T]}\|\rho(t)\|^{1/2}\Big(\|\zeta\|^{1/2}+\sup_{t\in[t_0,T]}\|l(t)\|^{1/2}\Big)
\\\notag&\qquad\qquad\qquad\qquad\qquad\qquad\qquad+\la\sup_{t\in[t_0,T]}\|w(t,\cdot)\|^{1/2}\Big(\|\rho_0\|^{1/2}+\sup_{t\in[t_0,T]}\|k(t)\|^{1/2}\Big)\Big).
\end{align}
Taking $\sup_{s\in(t_0,T]}\sup_{\{\xi:\|\xi\|=1\}}$ on both sides and rearranging, we obtain that 
\begin{align}\notag
\sup_{t\in[t_0,T]}\|\rho(t)\|&\le C\la \Big( \|\rho_0\|+\|\zeta\|+\sup_{t\in [t_0,T]}\left(\|k(t)\|+\|l(t)\|\right)\\\notag
&\qquad\qquad+\sup_{t\in[t_0,T]}\|w(t,\cdot)\|^{1/2}\Big(\|\rho_0\|^{1/2}+\sup_{t\in[t_0,T]}\|k(t)\|^{1/2}\Big).
\end{align}
From the last bound, \eqref{ASAF005}, and Young's inequality, we obtain that $\mathcal{A}$ is bounded. It also easily follows that any fixed point satisfies \eqref{ASAF004}. Since the system is linear, this bound also implies uniqueness.
%

\qed

Equipped with the last proposition, consider the solution of \eqref{1309} associated with $(\mu,u)$, the functions $k,l,\zeta\equiv 0$, and an initial condition $m_0\in\calM([d])$; denote it by $(m,v)$. That is, 
\begin{align}
\begin{cases}\label{1309b}
\frac{d}{dt} m_x(t)=\sum_{y\in[d]}m_y(t)a^*_x(y,\Delta_y u(t,\cdot))+\sum_{y\in[d]}\mu_y(t)\nabla_pa^*_x(y,\Delta_yu(t,\cdot))\cdot\Delta_yv(t,\cdot),\\
-\frac{d}{dt} v(t,x)=
\Delta_xv(t,\cdot)\cdot a^*(x,\Delta_xu(t,\cdot))
+\nabla_\eta f_2(x,\mu(t))\cdot m(t),\\
v(T,x)=\nabla_\eta g(x,\mu(T))\cdot m(T),\quad m(t_0)=m_0.
\end{cases}
\end{align}
From Proposition \ref{prop31}, $\sup_{t\in[t_0,T]}\|(m,v)\|$ is bounded by a constant that depends on $(\mu,u)$ and $m_0$, which in turn depends on the data of the problem given in Assumption \ref{assumption1300}.  
The next lemma provides a sensitivity result for $(m,v)$, similar to Lemma \ref{lem1301}.
\begin{lemma}\label{lem2301}
Fix  $0\le \hat t_0\le t_0\le T$ and $\mu_0,\hat\mu_0\in\calP([d])$. Let $(\mu,u)$ and $(\hat\mu,\hat u)$ be two solutions to \eqref{1300} with the initial conditions $\mu(t_0)=\mu_0$ and $\hat\mu(\hat t_0)=\hat\mu_0$. Also, let $(m,v)$ and $(\hat m,\hat v)$ be two solutions to \eqref{1309} associated with $(\mu,u)$ and $(\hat\mu,\hat u)$, respectively and satisfying the initial conditions $m(t_0)=\hat m(\hat t_0)=m_0$. 
Then, there is a constant $C>0$ independent of $t_0$, $\hat t_0$, $\mu_0$, and $\hat\mu_0$, such that 
\begin{align}\label{1300bbb}
&\sup_{t\in[t_0,T]}(\|m(t)-\hat m(t)\|
+\|v(t,\cdot)-\hat v(t,\cdot)\|
)
\le C\left(|t_0-\hat t_0|+\|\mu_0-\hat\mu_0\|\right),\\\label{ASAF001}
&\sup_{t\in[\hat t_0,t_0]}\|\hat m(t)-\hat m(t_0)\|\le  C|t_0-\hat t_0|.
\end{align}
\end{lemma}
{\bf Proof.}  Applying Proposition \ref{prop31} to $(\hat m, \hat v)$, we obtain the bound
\begin{align}\notag
\sup_{t\in[\hat t_0,T]}(\|\hat m(t)\|+\|\hat v(t,\cdot)\|)\le C,
\end{align} 
where hereafter in the rest of the proof, $C$ is a positive constant, independent of $t_0$, $\hat t_0$, $\mu_0$, and $\hat\mu_0$ (but may depend on $ m_0$), that can change from line to the next. 
Moreover, since $a^*$, $\mu$, $\nabla_p a^*$ are bounded, we obtain \eqref{ASAF001}.

Now, from \eqref{1300} and \eqref{1309b} it follows that $(\rho,w):=( m-\hat m,v-\hat v)$ satisfies \eqref{1309} on the time interval $[t_0,T]$ with the following data: $\rho_0=m_0-\hat m(t_0)$, 
\begin{align}\notag
k_x(t)&:=\sum_{y\in[d]}\mu_y(t)[\nabla_pa^*_x(y,\Delta_yu(t,\cdot))-\nabla_pa^*_x(y,\Delta_y\hat u(t,\cdot))]\cdot\Delta_x\hat v(t,\cdot)\\\notag
&\quad+\sum_{y\in[d]}[\mu_y(t)-\hat\mu_y(t)]\nabla_pa^*_x(y,\Delta_y\hat u(t,\cdot))\cdot\Delta_x\hat v(t,\cdot)\\\notag
&\quad+\sum_{y\in[d]}\hat m_y(t)[a^*_x(y,\Delta_y u(t,\cdot))-a^*_x(y,\Delta_y \hat u(t,\cdot))]\\\notag
l_x(t)&:=[a^*(x,\Delta_x u(t,\cdot))-a^*(x,\Delta_x \hat u(t,\cdot))]\cdot\Delta_x\hat v(t,\cdot)\\\notag
&\quad+[\nabla_\eta f_2(x,\mu(t))-\nabla_\eta f_2(x,\hat \mu(t))]\cdot\hat m(t),\\\notag
\zeta(x)&:=[\nabla_\eta g(x,\mu(T))-\nabla_\eta g(x,\hat\mu(T))]\cdot\hat m(T).
\end{align}
Again $\sum_{x\in[d]}a^*_x(y,p)=0$ implies that $\sum_{x\in[d]}k_x(t)=0$, hence the image of $k$ is in $\calM([d])$.
From Proposition \ref{prop31} it is sufficient to show that 
\begin{align}\notag
\|\rho_0\|+\|\zeta\|+\sup_{t\in [t_0,T]}\left(\|k(t)\|+\|l(t)\|\right)\le C(|t_0-\hat t_0|+\|\mu_0-\hat\mu_0\|).
 \end{align} 
The bound for $\|\rho_0\|$ follows from \eqref{ASAF001}. The rest of the bound follows from Lemma \ref{lem1301} observing that $\mu$, $\nabla_p a^*$, $\hat m$,and $\hat v$ are bounded, and $\nabla_p a^*$, $a^*$, $\nabla_\eta f_2$, and $\nabla_\eta g$ are Lipschitz continuous.

\qed

We now connect the functions $u$ and $v$.  
\begin{lemma}\label{lem1302}
Fix $t_0\in[0,T]$ and $\mu_0,\hat\mu_0\in\calP([d])$. Let $(u,\mu)$ and $(\hat u,\hat\mu)$ be solutions of \eqref{1300} with the initial conditions $\mu(t_0)=\mu_0$ and $\hat\mu(t_0)=\hat\mu_0$. Also, let $(v,m)$ be the solution of \eqref{1309} with the initial condition $m(t_0)=\hat\mu_0-\mu_0$. Then, there is a constant $C>0$ independent of $t_0$, $\mu_0$, and $\hat\mu_0$, such that 
\begin{align}\notag
\sup_{t\in[t_0,T]}&\left(\|\hat\mu(t)-\mu(t)-m(t)\|+\|\hat u(t,\cdot)-u(t,\cdot)-v(t,\cdot)\|
\right)
\le C\|\mu_0-\hat\mu_0\|^2.
\end{align}
\end{lemma}
{\bf Proof.} 
From \eqref{1300} and \eqref{1309b} it follows that $(\rho,w):=(\hat\mu-\mu-m,\hat u-u-v)$ satisfies \eqref{1309} with the following data: $\rho_0=0$, 
\begin{align}\notag
k_x(t)&:=\sum_{y\in[d]}(\hat \mu_y(t)-\mu_y(t))\nabla_pa^*_x(y,\Delta_yu(t,\cdot))\cdot(\Delta_y\hat u(t,\cdot)-\Delta_yu(t,\cdot))\\\notag
&\quad+\sum_{y\in[d]}\hat\mu_y(t)[a^*_x(y,\Delta_y\hat u(t,\cdot)-a^*_x(y,\Delta_y u(t,\cdot))\\\notag
&\qquad\qquad\qquad\qquad\qquad-\nabla_pa^*_x(y,\Delta_y u(t,\cdot)\cdot(\Delta_y\hat u(t,\cdot)-\Delta_y(t,\cdot))],\\\notag
l_x(t)&:=H_1(x,\Delta_x\hat u(t,\cdot))-H_1(x,\Delta_xu(t,\cdot))- a^*(x,\Delta_xu(t,\cdot))\cdot [\Delta_x\hat u(t,\cdot)-\Delta_x u(t,\cdot)]\\\notag
&\quad+f_2(x,\hat\mu(t))-f_2(x,\mu(t))-(\hat\mu(t)-\mu(t))\cdot\nabla_\eta f_2(x,\mu(t)),\\\notag
\zeta(x)&:=g(x,\hat\mu(T))-g(x,\mu(T))-(\hat\mu(T)-\mu(T))\cdot \nabla_\eta g(x,\hat\mu(T)).
\end{align}
Again the image of $k$ is in $\calM([d])$. 
The proof follows by Proposition \ref{prop31}, using the boundedness of $\nabla_p a^*$ and $\hat\mu$, the Lipschitz continuity of $\nabla_p a^*$, $H_1$, $f_2$, and $g$, and Lemma \ref{lem1301}. 

\qed

For every $(t_0,x,y,\mu_0)\in[0,T]\times[d]^2\times\calP([d])$ set $K_y(t_0,x,\mu_0):=v_{t_0,\mu_0,e_y}(t_0,x)$, where $v_{t_0,\mu_0,m}$ solves \eqref{1309} with the initial conditions $\mu(t_0)=\mu_0$ and $m(t_0)=m_0\in\R^d$. Denote $K=(K_y:y\in[d])$. From the linearity of the system \eqref{1309}, it follows that 
\begin{align}\label{1313}
v_{t_0,\mu_0,m_0}(t_0,x)=K(t_0,x,\mu_0)\cdot m_0.
\end{align}
Notice that from Lemma \ref{lem2301}, for every $0\le\hat t_0\le t_0\le T$ and $\mu_0,\hat\mu_0\in\calP([d])$, 
\begin{align}\notag
&|v_{t_0,\mu_0,e_y}(t_0,x)-v_{\hat t_0,\hat \mu_0,e_y}(\hat t_0,x)|\\\notag
&\quad
\le 
|v_{t_0,\mu_0,e_y}(t_0,x)-v_{\hat t_0,\hat \mu_0,e_y}( t_0,x)|
+
|v_{\hat t_0,\hat \mu_0,e_y}( t_0,x)-v_{\hat t_0,\hat \mu_0,e_y}(\hat t_0,x)|
\\\notag
&\quad\le C(|t_0-\hat t_0|+\|\mu_0-\hat\mu_0\|),
\end{align}
where $C>0$ is independent of $t_0$, $\hat t_0$, $\mu_0$, and $\hat\mu_0$.
Therefore, we get that 
\begin{align}\label{1313a}
\text{the mapping $\;(t_0,\mu_0)\mapsto K(t_0,x,\mu_0)\;$ is Lipschitz continuous.}
\end{align}
Also define the function $U:[0,T]\times[d]\times \calP([d])\to\R$ by 
$$U(t_0,x,\mu_0)=u_{t_0,\mu_0}(t_0,x),$$ 
where $u=u_{t_0,\mu_0}$ solves \eqref{1300} with the initial condition $\mu(t_0)=\mu_0$. We are now ready to prove the theorem. Namely, that $U$ uniquely solves \eqref{215} and has a Lipschitz $\eta$-gradient.

{\bf Proof.} [Proof of Theorem \ref{thm1300}] 
 {\bf Differentiability:} From Lemma \ref{lem1302} and \eqref{1313} we obtain that 
\begin{align}\label{1314}
\|U(t_0,\cdot,\hat\mu_0)-U(t_0,\cdot,\mu_0)-K(t_0,x,\mu_0)\cdot(\hat\mu_0-\mu_0)\|\le C\|\hat\mu_0-\mu_0\|^2,
\end{align}
and therefore, the function $U$ is differentiable w.r.t.~the measure component and $\nabla_\eta U(t_0,x,\mu_0)=K(t_0,x,\mu_0)$. Moreover, from \eqref{1313a}, we get the continuity of $\nabla_\eta U$. 

\skp {\bf Lipschitz continuity:} From \eqref{1313a} and \eqref{1314} the function $\nabla_\eta U$ is Lipschitz continuous w.r.t.~$\eta$, uniformly in $t$. That is, \eqref{225-ccc} holds.

\skp {\bf Satisfying \eqref{215}:} 
In this part we show that $U$ satisfies \eqref{215}. 
Clearly,
\begin{align}\label{1315}
\frac{1}{h}(U(t_0+h,x,\mu_0)-U(t_0,x,\mu_0))&=\frac{1}{h}(U(t_0+h,x,\mu_0)-U(t_0+h,x,\mu(t_0+h)))\\\notag
&\quad+\frac{1}{h}(U(t_0+h,x,\mu(t_0+h))-U(t_0,x,\mu_0)).
\end{align}
Set the matrix $\al\in\calQ^{d\times d}$ by $\al_{xy}(t,\eta):=a^*_y(x,\eta,\Delta_xU(s,\cdot,\eta))$. For every $s\in[0,1]$, denote by $\mu_s:=(1-s)\mu(t_0)+s\mu(t_0+h)$. Then,
\begin{align}\notag
&U(t_0+h,x,\mu(t_0+h))-U(t_0+h,x,\mu_0)\\\notag
&\qquad=\int_0^1\nabla_\eta U(t_0+h,x,\mu_s)\cdot(\mu(t_0+h)-\mu(t_0))ds\\\notag
&\qquad=\int_0^1\int_{t_0}^{t_0+h}\nabla_\eta U(t_0+h,x,\mu_s)\cdot(\al^\top(t,\mu(t))\mu(t))dtds,
\end{align}
where the matrix $\al^\top$ is the transpose of $\al$. For the last equality we used the differential equation for $\mu$ from \eqref{1300} and that $\Delta_xu(t,\cdot)=\Delta_xU(t,\cdot,\mu(t))$. 
Recall \eqref{1313a}, taking $\lim_{h\to 0}$ on both sides, we get that 
\begin{align}\notag
&\frac{1}{h}(U(t_0+h,x,\mu(t_0+h))-U(t_0+h,x,\mu_0))\\\notag
&\qquad=\nabla_\eta U(t_0,x,\mu_0)\cdot(\al^\top(t_0,\mu_0)\mu_0).
\end{align}
From the definition of $U$, we get
\begin{align}\notag
&\lim_{h\to 0}\frac{1}{h}(U(t_0+h,x,\mu(t_0+h))-U(t_0,x,\mu_0))\\\notag
&\quad=
\lim_{h\to 0}\frac{1}{h}(u(t_0+h,x)-u(t_0,x))=\frac{d}{dt} u(t_0,x)=H(x,\mu_0,\Delta_xU(t_0,\cdot,\mu_0)),
\end{align}
where in the last equality we used the differential equation for $u$ from \eqref{1300}, that $H=H_1+f_2$, and again the identity $\Delta_xu(t_0,\cdot)=\Delta_xU(t_0,\cdot,\mu_0)$. Combining the last two limits with \eqref{1315}, we get that $\frac{d}{dt} U(t_0,x,\mu_0)$ exists and that \eqref{215} holds.

\skp{\bf Uniqueness:} We now show that the master equation \eqref{215} has a unique solution. Let $\tilde U$ be a solution of the master equation, we show that $\tilde U=U$. Fix $\mu_0\in\calP([d])$ and let $\tilde \mu:[t_0,T]\to\calP([d])$ be the unique solution of 
\begin{align}\label{1319}
\begin{cases}
\frac{d}{dt}\tilde\mu^\top(t)=\tilde\mu^\top(t)\al(t,\tilde\mu(t)),\\
\tilde\mu(t_0)=\mu_0.
\end{cases}
\end{align}
Define the function $\tilde u:[t_0,T]\times[d]\to\R$ by $\tilde u(t,x)=\tilde U(t,x,\tilde\mu(t))$. Then,
using \eqref{1319} and the fact that $\tilde U$ solves \eqref{215}, we get that
\begin{align}\notag
\frac{d}{dt}\tilde u(t,x)=\partial_t\tilde U(t,x,\tilde\mu(t))+\nabla_\eta\tilde U(t,x,\tilde\mu(t))\cdot\partial_t\tilde\mu(t)=H_1(x,\Delta_x\tilde u(t,\cdot))+f_2(x,\tilde\mu(t)).
\end{align}
That is $(\tilde u,\tilde\mu)$ is a solution to \eqref{1300}, which by Lemma \ref{lem1301} has a unique solution. Therefore, $\tilde u=u$ and as a consequence $\tilde U =U$ and uniqueness of the solution of the master equation is established.

 \qed
 

\appendix
\section{Technical lemma} 
\beginsec

\begin{lemma}\label{lem_a1}
Let $D\subseteq\R^d$ and $f:[0,T]\times D\to\R$ whose gradient w.r.t.~$x\in D$ is Lipschitz-continuous uniformly in $t\in[0,T]$. Then there exists a function $g:[0,T]\times \R^{2d}\to\R$, bounded by the Lipchitz constant such that for every $x,y\in\R^d$,
\begin{align}\notag
f(t,x)=f(t,y)+(x-y)\cdot\nabla_x f(t,y)+g(t,x,y)\|x-y\|^2.
\end{align} 
\end{lemma}
{\bf Proof.} 
Fix $t\in[0,T]$ and $x\ne y$ and set
\begin{align}\label{ape1}
g(t,x,y)=\frac{f(t,x)-f(t,y)-(x-y)\cdot\nabla_x f(t,y)}{\|x-y\|^2}.
\end{align} 
For any $x,y\in\R^d$ set the function $h(t,\cdot)=h_{x,y}(t,\cdot):\R\to\R$, given by $h(t,u)=f(t,y+u(x-y)/\|x-y\|)$. By Lagrange's mean value theorem and since the function $f$ is differentiable, we get that there exists $u_{t,x,y}\in[0,\|x-y\|]$ such that 
\begin{align}\notag
f(t,x)-f(t,y)&=h(t,\|x-y\|)-h(t,0)=h'(u_{t,x,y})\|x-y\|\\\notag
&=(x-y)\cdot\nabla_x f(t,y+u_{t,x,y}(x-y)/\|x-y\|).
\end{align}
Plugging this into \ref{ape1}, we obtain
\begin{align}\notag
|g(t,x,y)|&\le\left\|\nabla_x f(t,y+u_{t,x,y}(x-y)/\|x-y\|)-\nabla_x f(t,y)\right\|/\|x-y\|\\\notag
&\le c u_{t,x,y}/\|x-y\|\le c,
\end{align} 
where $c$ is the Lipschitz constant of $\nabla_x f$.

\qed

\skp\noi{\bf Acknowledgement.} We are thankful to the anonymous AE and the two referees for their suggestions, which helped us to improve the presentation of the paper.

\footnotesize
\bibliographystyle{abbrv} 
\bibliography{bib_Asaf_IMS} 

\end{document}